\documentclass[12pt]{article}
\usepackage{latexsym}
    \title{{\bf Generalized Rationality and a
``Jacobi Identity'' for Intertwining Operator Algebras}}
    \author{Yi-Zhi Huang}
    \date{}
    \begin{document}
    \bibliographystyle{alpha}
    \maketitle

    \input amssym.def
    \input amssym
    \newtheorem{rema}{Remark}[section]
    \newtheorem{propo}[rema]{Proposition}
    \newtheorem{theo}[rema]{Theorem}
   \newtheorem{defi}[rema]{Definition}
    \newtheorem{lemma}[rema]{Lemma}
    \newtheorem{corol}[rema]{Corollary}
     \newtheorem{exam}[rema]{Example}
\newcommand{\binom}[2]{{{#1}\choose {#2}}}
	\newcommand{\nno}{\nonumber}
	\newcommand{\lbar}{\bigg\vert}
\newcommand{\mbar}{\mbox{\large $\vert$}}
	\newcommand{\p}{\partial}
	\newcommand{\dps}{\displaystyle}
	\newcommand{\bra}{\langle}
	\newcommand{\ket}{\rangle}
 \newcommand{\res}{\mbox{\rm Res}}
\renewcommand{\hom}{\mbox{\rm Hom}}
 \newcommand{\pf}{{\it Proof.}\hspace{2ex}}
 \newcommand{\epf}{\hfill$\Box$}
 \newcommand{\epfv}{\hfill$\Box$\vspace{1em}}
\newcommand{\nord}{\mbox{\scriptsize ${\circ\atop\circ}$}}
\newcommand{\wt}{\mbox{\rm wt}\ }
\newcommand{\clr}{\mbox{\rm clr}\ }

\begin{abstract}
We prove a generalized rationality property and  a new
identity that we call the ``Jacobi identity'' for intertwining operator
algebras.  Most of the main properties of genus-zero conformal field
theories, including the main properties of vertex operator algebras,
modules, intertwining operators, Verlinde algebras, and fusing and
braiding matrices, are incorporated into this identity. Together
with associativity and commutativity for
intertwining operators proved by
the author in \cite{H2} and \cite{H3}, 
the results of the present paper solve completely the
problem of finding a natural purely algebraic structure on the direct sum of
all inequivalent irreducible modules for a suitable vertex operator
algebra.  Two equivalent definitions of intertwining operator algebra
in terms of this Jacobi identity are given. 
\end{abstract}

\renewcommand{\theequation}{\thesection.\arabic{equation}}
\renewcommand{\therema}{\thesection.\arabic{rema}}
\setcounter{equation}{0}
\setcounter{rema}{0}
\setcounter{section}{-1}

\section{Introduction}

This paper is devoted to a study of certain general features of
representations of a vertex operator algebra. The basic definitions
that we shall need are recalled in Section 1 below. 

The direct sum of all inequivalent irreducible modules for a suitable
vertex operator algebra, equipped with intertwining operators, has a
natural algebraic structure called intertwining operator algebra (see
\cite{H2} and \cite{H4}).  Intertwining operator algebras are
multivalued analogues of vertex operator algebras. They were first
defined using the convergence property, associativity and
skew-symmetry as the main axiom. It is natural to ask whether
intertwining operator algebras have analogues of all the main
properties of vertex operator algebras and whether there is a
definition formulated using formal variables.

In the present paper, we prove an important property, a
generalized rationality property, for intertwining operator algebras,
and using this generalized rationality property and other properties
of intertwining operator algebras, we prove a new identity that we
call the ``Jacobi identity'' for intertwining operator algebras. This
Jacobi identity is a generalization of the ``Jacobi identity'' for
vertex operator algebras discovered by Frenkel, Lepowsky and Meurman
\cite{FLM} and independently by Borcherds (cf. \cite{FLM}), and can be
used as the main axiom for intertwining operator algebras.  (Frenkel,
Lepowsky and Meurman called their identity the ``Jacobi identity''
because of its subtle analogy with the classical Jacobi identity in
the definition of Lie algebra. We continue to use this terminology for
the same reason.)  Most of the main properties of genus-zero conformal
field theories, including the main properties of vertex operator
algebras, modules, intertwining operators, Verlinde algebras, and
fusing and braiding matrices, are incorporated into this identity.
Even in the special case that the intertwining operator algebra is a
vertex operator algebra, we obtain a new result that the rationality
property of a vertex operator algebra is a consequence of the other
``duality'' properties and the axioms for vertex operator algebras
excluding the Jacobi identity. In the course of formulating this new
result, we in particular review the elementary notions and basic
results of the theory of vertex operator algebras, making this paper
quite self-contained.


The notion of vertex algebra was introduced by Borcherds \cite{B} in a
mathematical setting, and a modified version, the notion of vertex
operator algebra, was introduced by Frenkel, Lepowsky and Meurman
\cite{FLM}.  In a physical context, Belavin, Polyakov and
Zamolodchikov formalized and studied 
the  operator product algebra structure in
conformal field theory (see \cite{BPZ}) and later physicists also
arrived at, though perhaps without a completely rigorous definition,
the notion of ``chiral algebra,'' a notion which essentially coincides
with the notion of vertex operator algebra (see e.g.  \cite{MS}).

In Borcherds' original definition, a vertex algebra is defined to be a
vector space equipped with infinitely many products, a vacuum vector
and an operator playing the role of the Virasoro algebra operator
$L(-1)$, satisfying suitable axioms. The main axiom in Borcherds'
original definition is an ``associator formula'' which specifies
precisely how the associativity of the infinitely many products
fails. In \cite{FLM}, Frenkel, Lepowsky and Meurman worked with the
generating function of the infinitely many products (the vertex
operator map) and discovered and exploited an identity for vertex
operators which they called the ``Jacobi identity'' or ``Jacobi-Cauchy
identity.'' In their definition of vertex operator algebra, they used
the Jacobi identity instead of the associator formula as the main
axiom. Borcherds' associator formula can in fact be recovered {}from
the Jacobi identity by taking the coefficient of $z_{1}^{-1}$ in this
identity, and it is an easy consequence of results obtained
in \cite{FHL} (see especially Remark 2.7.2 and Sections 3.2--3.4 of
\cite{FHL}; cf. \cite{H1} and \cite{L1}) 
that the associator formula also implies
the Jacobi identity when the skew-symmetry for vertex operators
\cite{B} and
the other axioms for vertex algebras or vertex operator algebras hold.

Motivated by families of concrete examples, the notion of vertex (operator)
algebra was generalized by Frenkel, Feingold and Ries \cite{FFR} and by
Dong and Lepowsky \cite{DL} to notions involving formal series with
nonintegral powers of the formal variables 
and one-dimensional representations of braid
groups. The most general such notion  is the notion of abelian
intertwining algebra (see \cite{DL1} and 
\cite{DL}). Again the main axiom in the
definition of abelian intertwining algebra is an identity, still called
the Jacobi identity, which is a generalization of the Jacobi identity
for vertex operator algebras.

These Jacobi identities play an important role in the study of vertex
operator algebras, abelian intertwining algebras and their
representations. In the presence of other axioms, they imply the
rationality property and the associativity, commutativity and
skew-symmetry relations, collectively called ``duality'' relations or
properties. These identities are canonical in a fundamental sense, and
numerous other properties of vertex operator algebras and abelian
intertwining algebras can be derived {}from them (see, for example,
\cite{FLM}, \cite{FHL} and \cite{DL}). In particular, most of the
algebraic relations that one uses, such as generalized commutator
relations, are special cases of the Jacobi identity, which is in fact
a generating function of an infinite family of identities.

These Jacobi identities are elementary identities involving only
formal series and linear maps. Though they might look unfamiliar to
people seeing them for the first time, anyone with a basic knowledge
of linear algebra and formal series will have no difficulty to begin
studying them.

In the representation theory of vertex operator algebras and algebraic
conformal field theory, intertwining operators (chiral vertex
operators) are one of the main
subjects of study. The important notions of fusion rule, fusing
matrix, braiding matrix, and Verlinde algebra for a vertex operator
algebra are all defined in terms of intertwining operators (see 
\cite{V}, \cite{TK}, \cite{MS} and \cite{FHL}).  In the
special case based on the $c=1/2$ Virasoro algebra minimal model,
Feingold, Ries and Weiner \cite{FRW} studied the intertwining
operators in detail and obtained some identities for these operators
by calculating the correlation functions directly. Their identities
are very different {}from the Jacobi identity for vertex operator
algebras or abelian intertwining operator algebras because the special
properties of the $c=1/2$ minimal Virasoro vertex operator algebra are
involved. For the same reason, it seems difficult to generalize their method
to study intertwining operators for other vertex operator algebras.

On the other hand, in \cite{H2} and \cite{H3}, the author showed using
the tensor product theory developed in \cite{HL1}--\cite{HL5},
\cite{H2} and \cite{HL8} that the intertwining operators for a
suitable vertex operator algebra satisfy generalizations of many
properties, including associativity and commutativity, of vertex
operators in vertex operator algebras.  Using these generalizations,
the author introduced the notion of intertwining operator algebra in
\cite{H4} and constructed rational genus-zero modular functors and
genus-zero holomorphic weakly conformal field theories 
{}from intertwining operator algebras in \cite{H9}. Abelian
intertwining algebras satisfying  grading-restriction conditions
are simple examples of intertwining operator algebras.

In \cite{H4}, the main axioms in the definition of intertwining
operator algebra are the convergence property, associativity and
skew-symmetry for intertwining operators. Based on the result proved in
\cite{H3} that the commutativity for intertwining operators follows
{}from the convergence property, associativity and skew-symmetry, it is
also remarked in \cite{H4} that the skew-symmetry can be replaced by
commutativity.  It is natural to ask whether there is an analogue of
the rationality property, whether there is a formulation
using only formal variables rather than convergence, multivalued
analytic functions and their expansions, and in particular, whether
there is an identity for intertwining operator algebras which
can be used as the main axiom and such that the Jacobi identities for
vertex operator algebras and abelian intertwining operator algebras
are special cases of this new identity.

This paper introduces a generalized rationality property and an  identity
that we call  the ``Jacobi identity''
for intertwining operator algebras.  The purpose of this paper is to
prove these properties (see Theorems \ref{rat}, \ref{associa},
\ref{jacobi}, and formulas 
(\ref{1.0}) and (\ref{1.1})), and to give a
formal-variable definition of intertwining operator algebra in terms
of this Jacobi identity (see Definitions \ref{4-1} and \ref{4-2}). 
In a forthcoming paper
\cite{H6}, we shall study the Jacobi identity and intertwining
operator algebras in detail.

For vertex operator algebras and abelian intertwining operator
algebras, the rationality property is either taken to be an axiom 
or derived from the Jacobi identity. For these algebras,
the Jacobi identities were first discovered using examples
constructed {}from lattices, because in the lattice case, the vertex
operators or the intertwining operators have explicit expressions
(see \cite{FLM}, \cite{DL1} and \cite{DL}, cf. also \cite{DL0}). 
But for nonabelian intertwining operator
algebras,  we do not have such simple examples as in
the lattice case for which the
intertwining operators can be expressed explicitly.
It is also in general very hard to verify directly the generalized
rationality property. In the present 
paper, the
generalized rationality property and the 
Jacobi identity are obtained {}from the convergence property, 
associativity, commutativity and 
other axioms for intertwining operator 
algebras.  

We give two forms of this Jacobi identity, one abstract and one
explicit. The abstract form (\ref{1.0}) will often be useful in studying
the axiomatic properties of the Jacobi identity and intertwining
operator algebras, while the explicit form (\ref{1.1}) will often be useful
for concrete calculations.

The results of the present paper provide a powerful method for the
study of the problem of constructing vertex operator algebras, abelian
intertwining operator algebras or intertwining operator algebras
starting {}from given vertex operator algebras and their modules. The
construction of the moonshine module vertex operator algebra 
by Frenkel, Lepowsky and Meurman
\cite{FLM1} \cite{FLM} can in fact be viewed in retrospect 
as the answer to an
important special case of this problem.  This problem was also studied
in \cite{DGM}, 
\cite{H2.5}, \cite{DLM} and \cite{L} in different situations.  The
construction of an abelian intertwining operator algebra underlying
the moonshine module given in \cite{H2.5} can be viewed as an application of
our method to a special but important case. In particular, a new and 
conceptual 
construction of vertex operator algebra structure on the moonshine module
was obtained in \cite{H2.5}.

In \cite{Z}, Zhu studied 
the vector spaces spanned by the $q$-traces
of products of  vertex operators for the irreducible modules 
for a rational vertex operator
algebra satisfying a certain finiteness condition, and proved that
these vector spaces are modular invariant (see also \cite{DLM2}). 
One step in Zhu's proof of the modular invariance
is to express the $q$-traces of products
of $n+1$ vertex operators for the irreducible modules as linear combinations
of the $q$-traces of products
of $n$ such operators,  with the expansions 
of suitable analytic functions as coefficients. 
This step reduces the proof of the
modular invariance of
the space spanned by the $q$-traces of products
of $n+1$ vertex operators  to
the proof of the same
for $n$ vertex operators, and thus to the proof
of the
modular invariance of
the space spanned by the $q$-traces of vertex operators (not products).
This step uses mainly the Jacobi identity
for vertex operators.
The Jacobi identity obtained in this paper allows us to generalize
this  step in Zhu's proof to the $q$-traces
of products of intertwining operators, and in particular, it allows us to
reduce the proof of the modular invariance of the spaces spanned by
the $q$-traces
of products of intertwining operators to the proof of the modular invariance
of the space spanned by the $q$-traces
of intertwining operators (not products). 
Using the Jacobi identity, we can also generalize  to
intertwining operator algebras the theory of Zhu's associative algebras and
Zhu's proof of the modular invariance of the 
space spanned by the
$q$-traces of vertex operators for the irreducible modules.
Thus we shall be able to prove that the
$q$-traces
of products of $n$ intertwining operators for any nonnegative
integer $n$  span a modular invariant
vector space. These modular invariant vector spaces are in fact
the fibers of the holomorphic vector bundles in the
genus-one modular functor for an associated
but yet-to-be-constructed conformal field
theory. These matters will be discussed in a future paper.

For the vertex operator algebras associated to the Virasoro algebra
and  affine Lie algebras, it was shown in \cite{H3} and \cite{HL7} that 
the convergence property, associativity and 
commutativity are satisfied by intertwining operators. 
Thus the intertwining operators for
these vertex operator algebras satisfy the Jacobi identity obtained in
this paper. We expect that this Jacobi identity will have many
applications to the representation theories of affine Lie algebras and
the Virasoro algebra. 

The Jacobi identity for vertex (operator) algebras has a reformulation
in terms of algebraic ${\cal D}$-modules (see \cite{BD} and
\cite{HL5}). It is an interesting  problem to 
reformulate the Jacobi identity for
intertwining operator algebras given in this paper using 
algebraic ${\cal D}$-modules or their generalizations. 

We recall basic notions and results in the representation theory 
of vertex operator algebras in Section 1  
(see \cite{FLM} and \cite{FHL} for more details). We also 
formulate our new result on the rationality property 
for vertex operator algebras in this section. In 
Section 2, we recall the definition and construction 
of intertwining operator algebras
and some basic results (see \cite{H2}, \cite{H3} and \cite{H4} 
for more details). 
In Section 3,  we introduce concepts and notations needed and
describe the Jacobi identity. We establish the generalized rationality
property and reformulate the duality properties in terms of formal
variables in Section 4. The Jacobi identity for
intertwining operator algebras follows immediately from these
generalized rationality and duality properties. 
Two
equivalent definitions of intertwining operator algebra 
in terms of the Jacobi
identity are given in Section 5.

\paragraph{Acknowledgment}  
I am indebted to David Kazhdan for noticing 
an inconsistency between 
the statements in a previous version of the present paper and the 
known expressions of certain correlation functions for WNZW models 
for ${\frak s}{\frak l}(2, {\Bbb C})$. This inconsistency 
led me to find and correct a mistake in the previous version.
I am grateful to Jim Lepowsky for comments.
This research is supported in part by an NSF grant 
DMS-9622961. 

\renewcommand{\theequation}{\thesection.\arabic{equation}}
\renewcommand{\therema}{\thesection.\arabic{rema}}
\setcounter{equation}{0}
\setcounter{rema}{0}

\section{Vertex operator algebras, modules and intertwining operators}

We recall basic notions and results in the theory of vertex operator
algebras and their representations in the present section; see
\cite{FLM}, \cite{FHL} and \cite{HL3}. In most of the section, we
follow \cite{FLM} and \cite{FHL} for these basic notions and
results. We also state a new result, proved later in this paper, on
rationality for vertex operator algebras (see Proposition
\ref{c&m->r}).

In this paper, all the variables $x$, $x_{0}, \dots$ are independent
commuting formal variables, and
for a vector space $W$ and a formal variable $x$, we
shall  denote the space of formal Laurent series $\sum_{n\in
{\Bbb Z}}w_{n}x^{n}$ ($w_{n}\in W$) in $x$ containing only finitely
many negative powers of $x$ by $W((x))$
and the space of formal sums of the form
$\sum_{n\in {\Bbb C}}w_{n}x^{n}$ ($w_{n}\in W$) by 
$W\{x\}$. We shall also use similar notations for series
with more than one formal variable. 
For any $f(x)\in W\{x\}$, we shall use $\res_{x}f(x)$ to denote the 
coefficient of $x^{-1}$ in $f(x)$. We shall also use the symbols $z,
z_{0}, \dots,$ which will denote complex numbers, {\it not} formal
variables.

For any $z\in {\Bbb C}$, we shall always choose $\log z$ so
that
$$
\log z=\log |z|+i\arg z\;\;\mbox{\rm with}\;\;0\le\arg z<2\pi.
$$
Given two multivalued functions $f_{1}$ and $f_{2}$ on a region, we
say that {\it $f_{1}$ and $f_{2}$ are equal} if on each simply
connected open subset of the region, for any single-valued branch of
$f_{1}$, there exists a single-valued branch of $f_{2}$ equal to it.

We use the formal expansion 
$$
\delta(x)=\sum_{n\in
{\Bbb Z}}x^{n}.
$$
This ``formal $\delta$-function'' has the following simple and
fundamental property:
For any $f(x)\in {\Bbb C}[x, x^{-1}]$,
$$
f(x)\delta(x)=f(1)\delta(x).
$$
This property has many important variants. For example, for any 
$$X(x_{1},
x_{2})\in (\mbox{End }W)[[x_{1}, x_{1}^{-1}, x_{2}, x_{2}^{-1}]]$$
(where $W$ is a vector space) such that 
\begin{equation}\label{0.1}
\lim_{x_{1}\to x_{2}}X(x_{1}, x_{2})=X(x_{1},
x_{2})\lbar_{x_{1}=x_{2}}=X(x_{2},
x_{2})
\end{equation}
exists, we have
$$
X(x_{1}, x_{2})\delta\left(\frac{x_{1}}{x_{2}}\right)=X(x_{2}, x_{2})
\delta\left(\frac{x_{1}}{x_{2}}\right).
$$
The existence of the ``algebraic limit'' defined in (\ref{0.1}) means that
for an arbitrary vector $w\in W$, the coefficient of each power of
$x_{2}$ in the formal expansion $X(x_{1}, x_{2})w\mbar_{x_{1}=x_{2}}$
is a finite sum.  We use the convention that negative powers of a
binomial are to be expanded in nonnegative powers of the second
summand. For example,
$$
x_{0}^{-1}\delta\left(\frac{x_{1}-x_{2}}{x_{0}}\right)=\sum_{n\in {\Bbb Z}}
\frac{(x_{1}-x_{2})^{n}}{x_{0}^{n+1}}=\sum_{m\in {\Bbb N},\; n\in {\Bbb Z}}
(-1)^{m}{{n}\choose {m}} x_{0}^{-n-1}x_{1}^{n-m}x_{2}^{m}.
$$
We have the following identities:
\begin{eqnarray*}
&{\dps x_{1}^{-1}\delta\left(\frac{x_{2}+x_{0}}{x_{1}}\right)=x_{2}^{-1}\left(
\frac{x_{1}-x_{0}}{x_{2}}\right),}&\\
&{\dps x_{0}^{-1}\delta\left(\frac{x_{1}-x_{2}}{x_{0}}\right)-
x_{0}^{-1}\delta\left(\frac{x_{2}-x_{1}}{-x_{0}}\right)=
x_{2}^{-1}\delta\left(\frac{x_{1}-x_{0}}{x_{2}}\right).}&
\end{eqnarray*}

We now recall the definition and basic ``duality''
properties of vertex operator algebras {from} \cite{FLM} or \cite{FHL}:

\begin{defi}
{\rm A {\it vertex operator algebra $($over ${\Bbb C}$$)$ of central
charge $c\in {\Bbb C}$} consists of the following data:

\begin{enumerate}

\item A ${\Bbb
Z}$-graded vector space 
$$
V=\coprod_{n\in {\Bbb Z}}V_{(n)}.
$$
For $v\in V_{(n)}$, $n$ is called the {\it weight} of $v$ and is
denoted by $\mbox{\rm wt}\ v$.

\item A linear map 
\begin{eqnarray*}
V&\to& (\mbox{\rm End}\; V)[[x, x^{-1}]]\nno\\
v\in V&\mapsto &Y(v, x)\in (\mbox{\rm End}\; V)[[x, x^{-1}]].
\end{eqnarray*}   

\item Two distinguished homogeneous vectors ${\bf 1}\in
V_{(0)}$ (the {\it vacuum}) and $\omega \in V_{(2)}$ (the {\it
Virasoro element}). 

\end{enumerate}

\noindent These data satisfy the following
conditions for $u, v \in V$: 

\begin{enumerate}

\item The {\it grading-restriction conditions}:
$$
\mbox{\rm dim }V_{(n)}<\infty\;\;\mbox{\rm for}\; n \in {\Bbb Z},
$$
$$
V_{(n)}=0\;\;\mbox{\rm for} \;n\; \mbox{\rm sufficiently small},
$$

\item The {\it lower truncation
condition}\footnote{The operator $Y_{n}(u)$ in this axiom 
is denoted $u_{n}$ in most of the 
literature on vertex operator algebras, including \cite{B} and
\cite{FLM}, where the notions of vertex algebra and vertex operator
algebra were first introduced.}: Let $Y_{n}(u)=\res_{x}x^{n}Y(u, x)\in
\mbox{\rm End}\;
V$, 
$n\in {\Bbb Z}$, that is, $Y(v, x)=\sum_{n\in {\Bbb
Z}}Y_{n}(u)x^{-n-1}$. 
Then 
$$
Y_{n}(u)v=0\;\;\mbox{\rm for}\;n\; \mbox{\rm sufficiently large}
$$
(or equivalently, $Y(u, x)v\in V((x))$).

\item The axioms for the vacuum:

\begin{enumerate}

\item The {\it identity property}
$$
Y({\bf 1}, x)=1\;\; (1\;\mbox{\rm on the right being the identity
operator}).
$$

\item The {\it creation property}:
$$
Y(v, x){\bf 1} \in V[[x]]\;\;\mbox{\rm and}\;\;\lim_{x\rightarrow
0}Y(v, x){\bf 1}=v
$$
(that is, $Y(v, x){\bf 1}$ involves only nonnegative integral powers
of $x$ and the constant term is $v$).

\end{enumerate}

\item The {\it Jacobi identity} (the
main axiom):
\begin{eqnarray*}
&x_{0}^{-1}\delta
\left({\displaystyle\frac{x_{1}-x_{2}}{x_{0}}}\right)Y(u, x_{1})Y(v,
x_{2})-x_{0}^{-1} \delta
\left({\displaystyle\frac{x_{2}-x_{1}}{-x_{0}}}\right)Y(v, x_{2})Y(u,
x_{1})&\nonumber \\ &=x_{2}^{-1} \delta
\left({\displaystyle\frac{x_{1}-x_{0}}{x_{2}}}\right)Y(Y(u, x_{0})v,
x_{2})&
\end{eqnarray*}
(note that when each expression in this identity 
 is applied to any element of
$V$, the coefficient of each monomial in the formal variables is a
finite sum; on the right-hand side, the notation $Y(\cdot, x_{2})$ is
understood to be extended in the obvious way to $V[[x_{0},
x^{-1}_{0}]]$).

\item The axioms for the Virasoro element: Let 
$$
L(n)=\res_{x}x^{n+1}Y(\omega, x)\;\; \mbox{\rm for} \;n\in{\Bbb Z}, 
\;\;{\rm
i.e.},\;\;Y(\omega, x)=\sum_{n\in{\Bbb Z}}L(n)x^{-n-2}.
$$

\begin{enumerate}

\item The {\it Virasoro algebra relations}: 
$$
[L(m), L(n)]=(m-n)L(m+n)+{\displaystyle\frac{1}{12}}
(m^{3}-m)\delta_{n+m,0}c
$$
for $m, n \in {\Bbb Z}$.

\item The  {\it $L(0)$-grading property}:
$$L(0)v=nv=(\mbox{\rm wt}\ v)v\;\;\mbox{\rm for}\;n \in {\Bbb
Z}\;\mbox{\rm and}\;v\in V_{(n)}.$$

\item The {\it  $L(-1)$-derivative property}:
$$\frac{d}{dx}Y(v,
x)=Y(L(-1)v, x).$$

\end{enumerate}
\end{enumerate}}
\end{defi}

We call $Y(v, x)$
the {\it vertex operator associated with} $v$.
The vertex operator algebra just defined is denoted by $(V, Y, {\bf 1},
\omega)$ (or simply by $V$).

Vertex operator algebras have important ``rationality,''
``commutativity'' and ``associativity'' properties, collectively
called ``duality'' properties.  These properties can in fact be used
as axioms replacing the Jacobi identity in the definition of vertex
operator algebra, as we now recall.

In the propositions below, ${\Bbb C}[x_{1}, x_{2}]_{S}$ is the ring of
rational functions obtained by inverting (localizing with respect to)
the products of (zero or more) elements of the set $S$ of nonzero
homogeneous linear polynomials in $x_{1}$ and $x_{2}$. Also,
$\iota_{12}$ (which might also be written as $\iota_{x_{1}x_{2}}$) is
the operation of expanding an element of ${\Bbb C}[x_{1}, x_{2}]_{S}$,
that is, a polynomial in $x_{1}$ and $x_{2}$ divided by a product of
homogeneous linear polynomials in $x_{1}$ and $x_{2}$, as a formal series
containing at most finitely many negative powers of $x_{2}$ (using
binomial expansions for negative powers of linear polynomials
involving both $x_{1}$ and $x_{2}$); similarly for $\iota_{21}$, and so
on. (The distinction between rational functions and formal Laurent
series is crucial.)

For any ${\Bbb Z}$-graded, or more generally, ${\Bbb C}$-graded, vector space 
$W=\coprod W_{(n)}$, we use the notation 
$$
W'=\coprod W_{(n)}^{*}
$$
for its graded dual.

\begin{propo}
{\bf (a) (rationality of products)} For $v$, $v_{1}$, $v_{2}\in V$ and
$v'\in V'$, the formal series
$\left\langle v', Y(v_{1}, x_{1})Y(v_{2}, x_{2})v\right\rangle,$ which
involves only finitely many negative powers of $x_{2}$ and only
finitely many positive powers of $x_{1}$, lies in the image of the map
$\iota_{12}$:
$$
\left\langle v', Y(v_{1}, x_{1})Y(v_{2}, x_{2})v\right\rangle
=\iota_{12}f(x_{1}, x_{2}),
$$
where the (uniquely determined) element $f\in {\Bbb C}[x_{1},
x_{2}]_{S}$ is of the form
$$
f(x_{1}, x_{2})={\displaystyle \frac{g(x_{1},
x_{2})}{x_{1}^{r}x_{2}^{s}(x_{1}-x_{2})^{t}}}
$$
for some $g\in {\Bbb C}[x_{1}, x_{2}]$ and $r, s, t\in {\Bbb Z}$.

{\bf (b) (commutativity)} We also have 
$$
\left\langle v', Y(v_{2}, x_{2})Y(v_{1}, x_{1})v\right\rangle
=\iota_{21}f(x_{1}, x_{2}).
$$
\end{propo}

\begin{propo}
{\bf (a) (rationality of iterates)} For $v$, $v_{1}$, $v_{2}\in V$ and
$v'\in V'$, the formal series
$\left\langle v', Y(Y(v_{1}, x_{0})v_{2}, x_{2})v\right\rangle,$
which involves only finitely many negative powers of $x_{0}$ and only
finitely many positive powers of $x_{2}$, lies in the image of the map
$\iota_{20}$:
$$
\left\langle v', Y(Y(v_{1}, x_{0})v_{2},
x_{2})v\right\rangle=\iota_{20}h(x_{0}, x_{2}),
$$
where the (uniquely determined) element $h\in {\Bbb C}[x_{0},
x_{2}]_{S}$ is of the form
$$
h(x_{0}, x_{2})={\displaystyle \frac{k(x_{0},
x_{2})}{x_{0}^{r}x_{2}^{s}(x_{0}+x_{2})^{t}}}
$$
for some $k\in {\Bbb C}[x_{0}, x_{2}]$ and $r, s, t\in {\Bbb Z}$.

{\bf (b)} The formal series 
$\left\langle v', Y(v_{1}, x_{0}+x_{2})Y(v_{2}, x_{2})v\right\rangle,$
which involves only finitely many negative powers of $x_{2}$ and only
finitely many positive powers of $x_{0}$, lies in the image of
$\iota_{02}$, and in fact
$$
\left\langle v', Y(v_{1}, x_{0}+x_{2})Y(v_{2},
x_{2})v\right\rangle=\iota_{02}h(x_{0}, x_{2}).
$$
\end{propo}

\begin{propo}[associativity]
We have the following equality of rational functions:
$$
\iota_{12}^{-1}\left\langle v', Y(v_{1}, x_{1})Y(v_{2}, x_{2})v\right\rangle
=(\iota_{20}^{-1}\left\langle v', Y(Y(v_{1}, x_{0})v_{2},
x_{2})v\right\rangle)\lbar_{x_{0}=x_{1}-x_{2}}.
$$
\end{propo}

\begin{propo}
In the presence of the other axioms, the Jacobi identity follows {from}
the rationality of products and iterates, commutativity and
associativity. In particular, in the definition of vertex operator
algebra, the Jacobi identity may be replaced by these properties.
\end{propo}

Another important property of vertex operator algebras is the 
skew-symmetry \cite{B}  which follows easily from the Jacobi identity (see
\cite{FHL}): For any $u, v\in V$,
$$Y(u, x)v=e^{xL(-1)}Y(v, -x)u.$$

We have:

\begin{propo}
In the presence of rationality and
the axioms for vertex operator algebras except
the Jacobi identity, 
any two of the three properties, associativity, commutativity and 
skew-symmetry, imply the other. In particular,
in the definition of vertex operator
algebra, the Jacobi identity may be replaced by the 
rationality, and any two of these three properties.
\end{propo}
\pf
We assume that rationality and
the axioms for vertex operator algebras except for
the Jacobi identity hold.
The skew-symmetry follows from the associativity and 
commutativity because these two properties imply the Jacobi identity and 
the Jacobi identity implies skew-symmetry. The associativity
being a consequence of
skew-symmetry and commutativity was proved in \cite{FHL}.
A similar argument proves that the  commutativity 
follows from skew-symmetry and associativity (see \cite{L1}). 
\epfv

We now give our new result for vertex operator algebras. 
A natural question is whether rationality follows from the other properties.
To answer and even to formulate this question precisely, we first have to 
formulate  associativity and commutativity without using 
rationality. 

\begin{description}

\item[Commutativity:]
For $v$, $v_{1}$, $v_{2}\in V$ and
$v'\in V'$, the series 
$$\left\langle v', Y(v_{1}, z_{1})Y(v_{2}, z_{2})v\right\rangle$$
and 
$$\left\langle v', Y(v_{2}, z_{2})Y(v_{1}, z_{1})v\right\rangle$$
are absolutely convergent in the regions $|z_{1}|>|z_{2}|>0$ and
$|z_{2}|>|z_{1}|>0$, respectively,
and their sums are analytic extensions of each other, that is, 
there exists a analytic function defined on a connected region 
containing both of the regions $|z_{1}|>|z_{2}|>0$ and
$|z_{2}|>|z_{1}|>0$ such that the sums of the two series above 
are the expansions of this function in the two regions.

\item[Associativity:]
For $v$, $v_{1}$, $v_{2}\in V$ and
$v'\in V'$, the series 
$$\left\langle v', Y(v_{1}, z_{1})Y(v_{2}, z_{2})v\right\rangle$$
and 
$$\left\langle v', Y(Y(v_{1}, z_{1}-z_{2})v_{2},
z_{2})v\right\rangle$$
are absolutely convergent in the regions $|z_{1}|>|z_{2}|>0$ and
$|z_{2}|>|z_{1}-z_{2}|>0$, respectively,
and their sums are analytic extensions of each other (or equivalently,
their sums are equal in the region $|z_{1}|>|z_{2}|>|z_{1}-z_{2}|>0$).
\end{description}

We can also formulate the rationality using complex variables:

\begin{description}

\item[Rationality of products:]
For $v$, $v_{1}$, $v_{2}\in V$ and
$v'\in V'$, the series 
$$\left\langle v', Y(v_{1}, z_{1})Y(v_{2}, z_{2})v\right\rangle$$
is the expansion in the region $|z_{1}|>|z_{2}|>0$ 
of a rational function in $z_{1}$ and $z_{2}$ 
with the only possible poles $z_{1}, z_{2}=0, \infty$ and $z_{1}=z_{2}$.

\item[Rationality of iterates:]
For $v$, $v_{1}$, $v_{2}\in V$ and
$v'\in V'$, the series 
$$\left\langle v', Y(Y(v_{1}, z_{1}-z_{2})v_{2},
z_{2})v\right\rangle$$
is the expansion in the region $|z_{2}|>|z_{1}-z_{2}|>0$
of a rational function in $z_{1}$ and $z_{2}$ 
with the only possible poles $z_{1}, z_{2}=0, \infty$ and $z_{1}=z_{2}$.
\end{description}

The following new result gives an answer to the question asked above:

\begin{propo}\label{c&m->r}
In the presence of the axioms for vertex operator algebras 
except for the Jacobi identity, both the rationality of products
and the rationality of iterates follow from  commutativity, associativity
and absolute convergence of the series
\begin{equation}\label{conv-3y}
\langle v', Y(v_{1}, z_{1})
Y(v_{2}, z_{2}) e^{z_{3}L(-1)}v\rangle
\end{equation}
in the region $|z_{1}|>|z_{2}|>|z_{3}|>0$ for $v_{1}, v_{2}, v\in V$,
$v'\in V'$.
\end{propo}

We shall prove this result in Section 4 as a consequence of 
generalized rationality of intertwining operator algebras.

\begin{rema}
{\rm Proposition \ref{c&m->r} above is
useful since in many cases it might not be 
easy to verify rationality.
One application of this result is to the geometric (or operadic)
formulation of the notion of vertex operator algebra given in 
\cite{H-1}, \cite{H0}, \cite{HL-1}, \cite{HL0} and \cite{H5}:
It can be shown easily using Proposition \ref{c&m->r} that 
in the  definition of geometric vertex operator  algebra (see 
\cite{H5}, Section 5.2), 
the meromorphicity axiom can be replaced by a 
much weaker meromorphicity axiom, assuming only the meromorphicity of
the maps corresponding to spheres with three punctures and local coordinates
vanishing at the punctures.}
\end{rema}

We have the following notion of module for 
vertex operator algebras:

\begin{defi} 
{\rm  Given a vertex operator algebra  
$(V,Y,{\bf 1},\omega )$,  a {\it module for  $V$}  
(or  $V$-{\it module}  or {\it representation space}) consists of
the following data:

\begin{enumerate}

\item A  ${\Bbb C}$-graded  vector space (graded by {\it weights})
$$
W =\coprod_{n\in {\Bbb C}}W_{(n)}.
$$
For $w \in  W_{(n)}$, $n$ is called the {\it weight} of $w$ and is
denoted by $\mbox{wt}\ w$.

\item A linear map  
\begin{eqnarray*}
V &\rightarrow & (\mbox{End}\ W)[[x,x^{-1}]]\nno\\
v\in V&\mapsto & Y(v, x)\in (\mbox{End}\ W)[[x,x^{-1}]].
\end{eqnarray*}

\end{enumerate}

\noindent These data 
satisfy ``all the defining 
properties of a vertex operator algebra that make sense,"  that is,
the following axioms for  
$u,v \in  V$  and  $w \in  W$,

\begin{enumerate}

\item The {\it grading-restriction conditions}:
$$
\dim \ W_{(n)} < \infty \;\;\mbox{for}\;\;n \in  {\Bbb C},
$$
$$
W_{(n)} = 0 \;\; \mbox{for}\;\;  n \;\; \mbox{whose real part is 
sufficiently small.}
$$

\item The {\it lower-truncation condition}\footnote{Similar
to the case of vertex operator algebras, the operator
$Y_{n}(v)$ in this axiom is usually denoted by $v_{n}$.}: Let 
$Y_{n}(v)=\res_{x}x^{n}Y(v, x)\in \mbox{\rm End}\;
W$, $n\in
{\Bbb Z}$, that is, $Y(v, x)=\sum_{n\in {\Bbb
Z}}(Y_{W})_{n}(v)x^{-n-1}$ 
(note that the sum is over  ${\Bbb Z}$,  not  ${\Bbb C}$).
Then 
$$
Y_n(v)w = 0 \;\; \mbox{for} \;\;n \;\; \mbox{sufficiently large}.
$$

\item The {\it identity property}:
$$
Y(\mbox{\bf 1},z) = 1.
$$

\item The {\it Jacobi identity}: 
\begin{eqnarray*}
&{\dps x^{-1}_0\delta \left( {x_1-x_2\over x_0}\right) Y(u,x_1)Y(v,x_2) - 
x^{-1}_0\delta \left( {x_2-x_1\over -x_0}\right) Y(v,x_2)Y(u,x_1)}&\nno\\
&{\dps = x^{-1}_2\delta \left( {x_1-x_0\over x_2}\right) Y(Y(u,x_0)v,x_2)}
\end{eqnarray*}
(note that on the right-hand side,
$Y(u,x_0)$  is the operator associated with  $V$).

\item The axioms for the vertex operator associated to the Virasoro
element in $V$: Let 
$$
L(n) = \res_{x}x^{n+1}Y(\omega, x)
\;\; \mbox{for}\;\;  n \in  {\Bbb Z}, \;\; 
\mbox{i.e.}, \;\; Y(\omega, x) =\sum_{n\in {\Bbb Z}}L(n)x^{-n-2}.
$$

\begin{enumerate}

\item The {\it Virasoro algebra 
relations}: 
$$
[L(m),L(n)] = (m-n)L(m+n) + {1\over 12}(m^3-m)\delta _{m+n,0}c
$$
for  $m,n \in  {\Bbb Z}$,  where $c$ is the central charge of $V$.

\item The {\it $L(0)$-grading property}:
$$
L(0)w = nw = (\mbox{wt}\ w)w \;\; \mbox{for}\;\; n \in  {\Bbb C}  
\;\;\mbox{and}\;\;  w \in  W_{(n)}.
$$

\item The {\it $L(-1)$-derivative property}:
$$
{d\over dx}Y(v,x) = Y(L(-1)v,x),
$$
where  $L(-1)$  is the operator on  $V$.
\end{enumerate}
\end{enumerate}}
\end{defi}

This completes the definition of module. We call $Y(v, x)$
the {\it vertex operator associated with} $v$.
We may denote the module just defined by
$(W,Y)$ (or simply by $W$). 
If necessary, we shall use $Y_{W}$ or similar notation to indicate
that the vertex operators concerned act on $W$. 

The following notion of   intertwining operator is introduced in \cite{FHL}:

\begin{defi} 
{\rm Let  $V$  be a vertex operator algebra and let  
$(W_1,Y_1)$,  $(W_2,Y_2)$  and  $(W_3,Y_3)$  be three  $V$-modules  (not 
necessarily distinct, and possibly equal to  $V)$.  An {\it intertwining 
operator of  type  
${W_3}\choose {W_1\ W_2} $}  is a linear map  
\begin{eqnarray*}
W_1 &\rightarrow & (\mbox{Hom}(W_2,W_3))\{x\}\nno\\
w \in W_{1}& \mapsto & {\cal Y}(w,x) \in (\mbox{Hom}(W_2,W_3))\{x\}
\end{eqnarray*}
such that ``all the defining properties of a module action that make sense 
hold."  That is, for  $v \in  V$,  $w_{(1)} \in  W_1$ and  
$w_{(2)} \in  W_2,$ ${\cal Y}$ satisfies the following axioms:

\begin{enumerate}

\item The {\it lower truncation condition}\footnote{The operator 
${\cal Y}_{n}(w)$ in this axiom is denoted $w_{n}$ in
\cite{FHL}, where the notion of intertwining operator was first
introduced in the representation theory of vertex operator algebras.}: Let
${\cal Y}_{n}(w)=\res_{x}x^{n+1}{\cal Y}(w, x)\in \mbox{\rm 
Hom}(W_2,W_3)$, $n\in {\Bbb C}$, that is, ${\cal Y}(w, x)=\sum_{n\in
{\Bbb C}}{\cal Y}_n(w)x^{-n-1}$. Then
$$
(w_{(1)})_{n}w_{(2)} = 0\;\;  \mbox{for}\;\; n \;\; \mbox{whose real part is 
sufficiently large.}
$$

\item  The  {\it Jacobi identity}: For the operators  
$Y_{1}(v,\cdot )$, 
$Y_{2}(v,\cdot )$, $Y_{3}(v,\cdot )$ and  
${\cal Y}(\cdot, x_{2})$  acting on the element  $w_{(2)}$, 
\begin{eqnarray*}
\lefteqn{\dps x^{-1}_0\delta \left( {x_1-x_2\over x_0}\right) 
Y_3(v,x_1){\cal Y}(w_{(1)},x_2)w_{(2)}}\nno\\
&&\hspace{2em}- x^{-1}_0\delta \left( {x_2-x_1\over -x_0}\right) 
{\cal Y}(w_{(1)},x_2)Y_2(v,x_1)w_{(2)}\nno \\
&&{\dps = x^{-1}_2\delta \left( {x_1-x_0\over x_2}\right) 
{\cal Y}(Y_1(v,x_0)w_{(1)},x_2)
w_{(2)}}
\end{eqnarray*}
(note that the first term on the left-hand side is algebraically meaningful 
because of the lower-truncation 
condition, and the other terms are meaningful by
the 
usual 
properties of modules; also note that this Jacobi identity involves integral 
powers of  $x_0$ and  $x_1$ and complex powers of  $x_2$).

\item The {\it $L(-1)$-derivative property}:
$$
{d\over dx}{\cal Y}(w_{(1)},x) = {\cal Y}(L(-1)w_{(1)}, x),
$$
where  $L(-1)$  is the operator acting on  $W_{1}$. 
\end{enumerate}}
\end{defi}

The intertwining operators of the
same type ${W_{3}}\choose {W_{1}\ W_{2}}$ form a vector space, which we
denote  by ${\cal
V}^{W_{3}}_{W_{1}W_{2}}$. The dimension of this vector space 
is called the {\it fusion rule} for $W_{1}$,
$W_{2}$ and $W_{3}$ and is denoted by 
${\cal N}^{W_{3}}_{W_{1}W_{2}}$.

Let ${\cal Y}$ be an intertwining operator of type 
${W_{3}}\choose {W_{1}W_{2}}$. For any complex number $\zeta$ and any 
$w_{(1)}\in W_{1}$, 
${\cal Y}(w_{(1)}, y)\lbar_{y^{n}=e^{n\zeta}x^{n}, \ n\in {\Bbb C}}$ is 
a well-defined element of Hom$(W_{2}, W_{3})\{ x\}$. We denote
this element by ${\cal Y}(w_{(1)}, e^{\zeta}x)$. Note that this element 
depends on $\zeta$, not on just $e^{\zeta}$.
Given any $r\in {\Bbb Z}$, we define 
$\Omega_{r}({\cal Y}):W_2\otimes W_1 \rightarrow  W_3\{ x\}$
by the formula
$$
\Omega_{r}({\cal Y})(w_{(2)},x)w_{(1)} = e^{xL(-1)}
{\cal Y}(w_{(1)},e^{ (2r+1)\pi i}x)w_{(2)}.
$$
for  $w_{(1)}\in W_{1}$ and $w_{(2)}\in W_{2}$. The following result
is proved in \cite{HL3}:

\begin{propo}
The operator  $\Omega_{r}({\cal Y})$  is an intertwining 
operator of type  ${W_{3}}\choose {W_{2} \ W_{1}}$.  Moreover,
$$
\Omega_{-r-1}(\Omega_{r}({\cal Y}))=\Omega_{r}(\Omega_{-r-1}({\cal Y}))
 = {\cal Y}.
$$
In particular, the correspondence  ${\cal Y} \mapsto  
\Omega_{r}({\cal Y})$  defines a linear isomorphism {from}  
${\cal V}^{W_{3}}_{W_{1}W_{2}}$ to  
${\cal V}^{W_{3}}_{W_{2}W_{1}}$,  and we have
$$
{\cal N}^{W_{3}}_{W_{1}W_{2}} = {\cal N}^{W_{3}}_{W_{2}W_{1}}.
$$
\end{propo}

We shall also need the notion of rational vertex operator algebra
as defined in \cite{HL1} and \cite{HL2} 
in the next section. (There are different
definitions of ``rational vertex operator algebra'' used
for different purposes.)

\begin{defi}
{\rm A  vertex operator algebra $V$ is {\it rational} if it
satisfies the following conditions:
\begin{enumerate}
\item There are only finitely many irreducible $V$-modules (up to equivalence).
\item Every $V$-module is completely reducible (and is in particular a 
{\it finite} direct sum of irreducible modules).
\item All the fusion rules for $V$ are finite (for triples of irreducible
modules and hence arbitrary modules).
\end{enumerate}
}
\end{defi}

\renewcommand{\theequation}{\thesection.\arabic{equation}}
\renewcommand{\therema}{\thesection.\arabic{rema}}
\setcounter{equation}{0}
\setcounter{rema}{0}

\section{Intertwining operator algebras}

In \cite{H4}, the author introduced a notion of intertwining
operator algebra. In terms of vertex operator algebras, modules and
intertwining operators, the definition of intertwining operator
algebras can be reformulated as follows:

\begin{defi}\label{ioa}
{\rm An {\it intertwining operator algebra of central charge 
$c\in {\Bbb C}$} consists of the following data:

\begin{enumerate}

\item A vector space
$$W=\coprod_{a\in
{\cal A}}W^{a}$$
graded
by a finite set ${\cal A}$
containing a special element $e$
(graded  by {\it color}).

\item A vertex operator algebra
structure of central charge $c$ 
on $W^{e}$, and a $W^{e}$-module structure on $W^{a}$ for 
each $a\in {\cal A}$.

\item A subspace ${\cal V}_{a_{1}a_{2}}^{a_{3}}$ of 
the space of all intertwining operators of type 
${W^{a_{3}}\choose W^{a_{1}}W^{a_{2}}}$ for  each triple
$a_{1}, a_{2}, a_{3}\in {\cal A}$.

\end{enumerate}

\noindent These data satisfy the
following axioms for $a_{1}, a_{2}, 
a_{3}, a_{4}, a_{5}, 
a_{6}\in {\cal A}$, 
$w_{(a_{i})}\in W^{a_{i}}$, $i=1, 2, 3$, and $w_{(a_{4})}'
\in W'_{a_{4}}$:

\begin{enumerate}

\item The $W^{e}$-module structure on $W^{e}$ is the adjoint module
structure. For any $a\in {\cal A}$, the space ${\cal V}_{ea}^{a}$ is the
one-dimensional vector space spanned by the vertex operators for the
$W^{e}$-module $W^{a}$.
For any $a_{1}, a_{2}\in {\cal A}$ such that
$a_{1}\ne a_{2}$, ${\cal V}_{ea_{1}}^{a_{2}}=0$.

\item {\it Convergence properties}: for any $m\in {\Bbb Z}_{+}$,
$a_{i}, b_{i}, \mu_{i}\in {\cal A}$, $w_{(a_{i})}
\in W^{a_{i}}$, ${\cal Y}_{i}\in {\cal 
V}_{a_{i}\;b_{i+1}}^{\mu_{i}}$, $i=1, \dots, m$, $w_{(\mu_{1})}'
\in (W^{\mu_{1}})'$ and 
$w_{(b_{m})}\in W^{b_{m}}$, the series
\begin{equation}\label{conv-pr}
\langle w_{(\mu_{1})}', {\cal Y}_{1}(w_{(a_{1})}, x_{1})
\cdots{\cal Y}_{m}(w_{(a_{m})},
x_{m})w_{(b_{m})}\rangle_{W^{\mu_{1}}}\mbar_{x^{n}_{i}=e^{n\log z_{i}},
i=1, \dots, m, n\in {\Bbb R}}
\end{equation}
is absolutely convergent when $|z_{1}|>\cdots >|z_{m}|>0$, and for 
any ${\cal Y}_{1}\in {\cal 
V}_{a_{1}a_{2}}^{a_{5}}$ and ${\cal Y}_{2}\in {\cal 
V}_{a_{5}a_{3}}^{a_{4}}$, the series
\begin{equation}\label{conv-it}
\langle w_{(a_{4})}', {\cal Y}_{2}({\cal Y}_{1}(w_{(a_{1})}, 
x_{0})w_{(a_{2})},
x_{2})w_{(a_{3})}\rangle_{W^{a_{4}}}
\mbar_{x^{n}_{0}=e^{n\log (z_{1}-z_{2})},
x^{n}_{2}=e^{n\log z_{2}}, n\in {\Bbb R}}
\end{equation}
is absolutely convergent when
$|z_{2}|>|z_{1}-z_{2}|>0$.

\item {\it Associativity}: for any ${\cal Y}_{1}\in {\cal 
V}_{a_{1}a_{5}}^{a_{4}}$ and ${\cal 
V}_{a_{2}a_{3}}^{a_{5}}$, there exist ${\cal Y}^{a}_{3}
\in {\cal 
V}_{a_{1}a_{2}}^{a}$ and ${\cal Y}^{a}_{4}\in {\cal 
V}_{aa_{3}}^{a_{4}}$ for all $a\in {\cal A}$
such that the (multivalued) analytic function 
$$\langle w_{(a_{4})}', 
{\cal Y}_{1}(w_{(a_{1})}, x_{1}){\cal Y}_{2}(w_{(a_{2})}, 
x_{2})w_{(a_{3})}\rangle_{W}\mbar_{x_{1}=z_{1},
x_{2}=z_{2}}$$
defined in the region
$|z_{1}|>|z_{2}|>0$ 
and the (multivalued) analytic function
$$\sum_{a\in {\cal A}}\langle w_{(a_{4})}', {\cal Y}^{a}_{4}
({\cal Y}^{a}_{3}(w_{(a_{1})},
x_{0})w_{(a_{2})}, x_{2})w_{(a_{3})}\rangle_{W^{a_{4}}}
\mbar_{x_{0}=z_{1}-z_{2},
x_{2}=z_{2}}$$ 
defined in the region
$|z_{2}|>|z_{1}-z_{2}|>0$ are equal in the intersection
$|z_{1}|> |z_{2}|>|z_{1}-z_{2}|>0$.

\item {\it Skew-symmetry}: The restriction of  $\Omega_{-1}$ to
${\cal V}_{a_{1}a_{2}}^{a_{3}}$ is an isomorphism from
${\cal V}_{a_{1}a_{2}}^{a_{3}}$ to 
${\cal V}_{a_{2}a_{1}}^{a_{3}}$.

\end{enumerate}}
\end{defi}

To make our study slightly easier, we shall assume in the present
paper that intertwining operator algebras also satisfy the following
{\it additional grading-restriction condition}: 
For any $a\in {\cal A}$, there exists
$h_{a}\in {\Bbb R}$ such that $(W^{a})_{(n)}=0$ for $n\not
\in h_{a}+{\Bbb Z}$. But this is in fact a very minor 
restriction and in addition, all the results on intertwining operator
algebras satisfying this additional property can be generalized easily
to intertwining operator algebras not satisfying this condition.

The intertwining operator algebra just defined above is denoted by 
$$(W, 
{\cal A}, \{{\cal V}_{a_{1}a_{2}}^{a_{3}}\}, {\bf 1}, 
\omega)$$ or simply $W$. 

The following result is proved in \cite{H3}:

\begin{propo}
Let $(W, 
{\cal A}, \{{\cal V}_{a_{1}a_{2}}^{a_{3}}\}, {\bf 1}, \omega)$ 
be an intertwining operator algebra. Then we have the following
{\it commutativity}: for any ${\cal Y}_{1}\in {\cal 
V}_{a_{1}a_{5}}^{a_{4}}$ and ${\cal 
V}_{a_{2}a_{3}}^{a_{5}}$, there exist ${\cal Y}^{a}_{5}
\in {\cal 
V}_{a_{2}a}^{a_{4}}$ and ${\cal Y}^{a}_{6}\in {\cal 
V}_{a_{2}a_{3}}^{a}$ for all $a\in {\cal A}$
such that the (multivalued) analytic function 
$$\langle w_{(a_{4})}', 
{\cal Y}_{1}(w_{(a_{1})}, x_{1}){\cal Y}_{2}(w_{(a_{2})}, 
x_{2})w_{(a_{3})}\rangle_{W}\mbar_{x_{1}=z_{1},
x_{2}=z_{2}}$$
defined in the region
$|z_{1}|>|z_{2}|>0\}$ 
and the (multivalued) analytic function
$$\sum_{a\in {\cal A}}\langle w_{(a_{4})}', 
{\cal Y}_{5}(w_{(a_{2})}, x_{2}){\cal Y}_{6}(w_{(a_{1})}, 
x_{1})w_{(a_{3})}\rangle_{W^{a_{4}}}\mbar_{x_{1}=z_{1},
x_{2}=z_{2}}$$
defined in the region $|z_{2}|>|z_{1}|>0\}$ 
 are analytic extensions of each other.
In the definition of intertwining operator algebra, skew-symmetry 
can be replaced by  commutativity.
\end{propo}

\begin{rema}
{\rm In fact, the absolute convergence of (\ref{conv-it}) follows {}from 
the absolute convergence of (\ref{conv-pr}). See \cite{H2} for a proof.}
\end{rema}

The problem of constructing intertwining operator algebras from
modules and intertwining operators for 
suitable vertex operator algebras was solved implicitly in \cite{H2},
in which the associativity of
intertwining operators, which is the only nontrivial property,
was proved. To formulate this result precisely, we need 
the following notion introduced in \cite{H2}:

\begin{defi}
{\rm Let $V$ be a vertex operator algebra. We say that products of
intertwining operators
for $V$ satisfy the {\it
convergence and extension property} 
if for any intertwining operators
${\cal Y}_{1}$ and ${\cal Y}_{2}$ of types
${W_{0}}\choose {W_{1}W_{4}}$ and ${W_{4}}\choose
{W_{2}W_{3}}$, respectively,
there exists 
an integer $N$
(depending only on ${\cal Y}_{1}$ and ${\cal Y}_{2}$), and 
for any $w_{(1)}\in W_{1}$,
$w_{(2)}\in W_{2}$, $w_{(3)}\in W_{3}$, $w'_{(0)}\in W'_{0}$, there exist
$j\in {\Bbb N}$, $r_{i}, s_{i}\in {\Bbb R}$, $i=1, \dots, j$, and analytic 
functions $f_{i}(z)$ on $|z|<1$, $i=1, \dots, j$, 
satisfying
$$
\wt w_{(1)}+\wt w_{(2)}+s_{i}>N,\;\;\;i=1, \dots, j,
$$
such that
$$
\langle w'_{(0)}, {\cal Y}_{1}(w_{(1)}, x_{1})
{\cal Y}_{2}(w_{(2)}, x_{2})w_{(3)}\rangle
\lbar_{x_{1}= z_{1}, \;x_{2}=z_{2}}
$$
is absolutely 
convergent when $|z_{1}|>|z_{2}|>0$ and can be analytically extended to  
the multivalued analytic function
$$
\sum_{i=1}^{j}z_{2}^{r_{i}}(z_{1}-z_{2})^{s_{i}}
f_{i}\left(\frac{z_{1}-z_{2}}{z_{2}}\right)
$$
when $|z_{2}|>|z_{1}-z_{2}|>0$.}
\end{defi}

Note that if the associativity of intertwining operators holds for 
$V$ (see the associativity axiom in Definition \ref{ioa}), then products 
of intertwining operators satisfy the convergence and extension property,
so that this condition is necessary for the  
associativity of intertwining operators.

We shall also use the concept of generalized module, as defined in
\cite{HL2}.  Let $V$ be a vertex operator algebra. A {\it generalized
$V$-module} is a ${\Bbb C}$-graded vector space equipped with a vertex
operator map satisfying all the axioms for a $V$-module except for the two
grading-restriction conditions. If there exists $N\in {\Bbb Z}$ such that
the homogeneous subspace of weight $n$ of a generalized $V$-module is
$0$ when the real part of $n$ is less than $N$, the generalized
$V$-module is said to be {\it lower truncated}.

Now we have the following immediate consequence of the results
obtained in \cite{H2}:

\begin{theo}\label{cons-ioa}
Let $V$ be a rational vertex operator algebra and 
$${\cal A}=
\{W^{1}, \dots,
W^{m}\}$$ 
a set of representatives of equivalence classes of
irreducible $V$-modules. Assume that every irreducible $V$-module
is ${\Bbb R}$-graded, products of
intertwining operators
for $V$ satisfy the
convergence and extension property and every finitely-generated
lower-truncated generalized $V$-module is a $V$-module. Then the 
${\Bbb R}$-graded vector space
$\coprod_{i=1}^{m}W^{i}$ together with the finite set ${\cal A}$ with
the special element $e=V$, 
the spaces ${\cal V}_{a_{1}a_{2}}^{a_{3}}$, $a_{1}, a_{2}, a_{3}\in
{\cal A}$, the vacuum ${\bf 1}\in V$ and the Virasoro element
$\omega\in V$ is an intertwining operator algebra.
\end{theo}

In the case of familiar vertex operator algebras (e.g., the vertex
operator algebras associated to the Virasoro algebra and affine Lie
algebras), the conditions in this theorem were verified in \cite{H3}
and \cite{HL7} and thus we obtain intertwining operator algebras from
representations of these vertex operator algebras.

\renewcommand{\theequation}{\thesection.\arabic{equation}}
\renewcommand{\therema}{\thesection.\arabic{rema}}
\setcounter{equation}{0}
\setcounter{rema}{0}

\section{The Jacobi identity}

Let $(W, 
{\cal A}, \{{\cal V}_{a_{1}a_{2}}^{a_{3}}\}, \omega)$
be a vertex operator algebra. 
Since we assume that intertwining operator algebras satisfy 
the additional grading-restriction condition, we know that
there exist $h_{a}\in {\Bbb R}$, $a\in {\cal A}$,  
such that for any $a\in {\cal A}$,
$W_{a}=\coprod_{n\in {\Bbb Z}}(W_{a})_{(h_{a}+n)}$. 
For any $a_{1},
a_{2}, a_{3}\in {\cal A}$,  an
intertwining operator ${\cal Y}$ of type ${a_{3}\choose a_{1}a_{2}}$
is in fact a map {}from $W_{a_{1}}\otimes W_{a_{2}}$ to
$x^{h_{a_{3}}-h_{a_{1}}-h_{a_{2}}}W_{a_{3}}((x))$. See \cite{FHL} for
detailed discussions.  Let
$$H=\{h_{a}\;|\; a\in {\cal A}\}$$ 
and for any $a_{1}, a_{2}\in {\cal A}$, let
$${\Bbb P}(a_{1}, a_{2})=\{h_{a_{1}}+h_{a_{2}}-h+{\Bbb Z}\in {\Bbb
R}/{\Bbb Z}\;|\; h\in H\}.$$

The skew-symmetry isomorphism
$\Omega_{-1}(a_{1}, a_{2}; a_{3})$ 
for all $a_{1}, a_{2}, a_{3}\in {\cal A}$ give an isomorphism 
$$\Omega_{-1}: \coprod_{a_{1}, a_{2}, a_{3}\in {\cal A}}
{\cal V}_{a_{1}a_{2}}^{a_{3}}\to \coprod_{a_{1}, a_{2}, a_{3}\in {\cal A}}
{\cal V}_{a_{1}a_{2}}^{a_{3}}$$
and we still call this isomorphism the {\it skew-symmetry isomorphism}.
In this paper,  for simplicity, we shall omit subscript $-1$ in 
$\Omega_{-1}(a_{1}, a_{2}; a_{3})$,
$a_{1}, a_{2}, a_{3}\in {\cal A}$, and 
in $\Omega_{-1}$ and denote them 
simply by $\Omega(a_{1}, a_{2}; a_{3})$ and $\Omega$, 
respectively.

Note that in the associativity property, ${\cal Y}^{a}_{3}$ and
${\cal Y}^{a}_{4}$, $a\in {\cal A}$, are not unique. But 
if we require that 
$$\langle w'_{(a_{4})}, 
{\cal Y}_{1}(w_{(a_{1})}, x_{1}){\cal Y}_{2}(w_{(a_{2})}, 
x_{2})w_{(a_{3})}\rangle_{W_{a_{4}}}\mbar_{x_{1}^{n}=e^{n \log z_{1}},
x_{2}^{n}=e^{n\log z_{2}}}$$
is equal to 
$$\sum_{a\in {\cal A}}\langle w', {\cal Y}^{a}_{4}({\cal Y}^{a}_{3}(w_{1},
x_{0})w_{2}, x_{2})w_{3}\rangle_{W}\mbar_{x_{0}^{n}=e^{n \log (z_{1}-z_{2})},
x_{2}^{n}=e^{n\log z_{2}}}$$ 
when $z_{1}$ and $z_{2}$ are positive real numbers  satisfying
$z_{1}>z_{2}>z_{1}-z_{2}>0$, then ${\cal Y}^{a}_{3}\otimes 
{\cal Y}^{a}_{4}\in {\cal V}_{a_{1}a_{2}}^{a}\otimes {\cal
V}_{aa_{3}}^{a_{4}}$, $a\in {\cal A}$, are uniquely determined. 
Thus for any $a_{1}, a_{2}, a_{3}, a_{4}\in {\cal A}$, the associativity
of intertwining operators together with the above requirement 
gives a {\it fusing isomorphism} 
$${\cal F}(a_{1}, a_{2}, a_{3}; a_{4}):
\coprod_{a\in {\cal A}}{\cal V}_{a_{1}a}^{a_{4}}\otimes 
{\cal V}_{a_{2}a_{3}}^{a}
\to \coprod_{a\in {\cal A}}{\cal V}_{a_{1}a_{2}}^{a}\otimes {\cal
V}_{aa_{3}}^{a_{4}}.$$ 
The fusing isomorphisms for all $a_{1}, a_{2}, a_{3}, a_{4}\in {\cal A}$
give an isomorphism 
$${\cal F}:
\coprod_{a_{1}, a_{2}, a_{3}, a_{4},
a_{5}\in {\cal A}}{\cal V}_{a_{1}a_{5}}^{a_{4}}\otimes 
{\cal V}_{a_{2}a_{3}}^{a_{5}}
\to \coprod_{a_{1}, a_{2}, a_{3}, a_{4},
a_{5}\in {\cal A}}{\cal V}_{a_{1}a_{2}}^{a_{5}}\otimes {\cal
V}_{a_{5}a_{3}}^{a_{4}},$$
which will  still be called the {\it fusing isomorphism}.

For any $a_{1}, a_{2}, a_{3}, a_{4}\in {\cal A}$, 
we have a {\it braiding isomorphism}  
\begin{eqnarray}\label{braiding}
\lefteqn{{\cal B}(a_{1}, a_{2}; a_{3}, a_{4})}\nno\\
&&={\cal F}^{-1}(a_{2},
a_{1}, a_{3}; a_{4})\circ \left(\coprod_{a\in 
{\cal A}}(\Omega(a_{1}, a_{2}; a)
\otimes I_{{\cal V}_{aa_{3}}^{a_{4}}})\right)\circ 
{\cal F}(a_{1},
a_{2}, a_{3}; a_{4})\nno\\
&&
\end{eqnarray}
{}from $$\coprod_{a\in {\cal A}}{\cal
V}_{a_{1}a}^{a_{4}}\otimes {\cal V}_{a_{2}a_{3}}^{a}$$ to
$$\coprod_{a\in
{\cal A}}{\cal V}_{a_{2}a}^{a_{4}}\otimes {\cal V}_{a_{1}a_{3}}^{a}.$$
The braiding isomorphisms for all $a_{1}, a_{2}, a_{3}, a_{4}\in {\cal A}$
give an isomorphism
$${\cal B}: \coprod_{a_{1}, a_{2}, a_{3}, a_{4}, a_{5}\in {\cal A}}{\cal
V}_{a_{1}a_{5}}^{a_{4}}\otimes {\cal V}_{a_{2}a_{3}}^{a_{5}}
\to \coprod_{a_{1}, a_{2}, a_{3}, a_{4}, a_{5}\in {\cal A}}{\cal
V}_{a_{2}a_{5}}^{a_{4}}\otimes {\cal V}_{a_{1}a_{3}}^{a_{5}},$$
which will still be called the {\it braiding isomorphism}.

The skew-symmetry,  fusing and braiding isomorphisms induce 
isomorphisms between vector spaces containing the domains and images
of these isomorphisms. These induced isomorphisms are not
independent. To describe the relations satisfied by
these induced isomorphisms, we need to introduce notations for certain
particular induced isomorphisms. 

The vector spaces
$$\coprod_{a_{1}, a_{2}, a_{3}, a_{4}, a_{5}, a_{6}, a_{7}
\in {\cal A}}{\cal
V}_{a_{1}a_{6}}^{a_{5}}\otimes {\cal V}_{a_{2}a_{7}}^{a_{6}}
\otimes {\cal V}_{a_{3}a_{4}}^{a_{7}}$$ 
and 
$$\coprod_{a_{1}, a_{2}, a_{3}, a_{4}, a_{5}, a_{6}, a_{7}
\in {\cal A}}{\cal
V}_{a_{1}a_{2}}^{a_{6}}\otimes {\cal V}_{a_{6}a_{7}}^{a_{5}}
\otimes {\cal V}_{a_{3}a_{4}}^{a_{7}}$$ 
are isomorphic to 
$$\coprod_{a_{1}, a_{2}, a_{3}, a_{4}, a_{5}, a_{7}
\in {\cal A}}\left(\coprod_{a_{6}\in {\cal A}}{\cal
V}_{a_{1}a_{6}}^{a_{5}}\otimes {\cal V}_{a_{2}a_{7}}^{a_{6}}\right)
\otimes {\cal V}_{a_{3}a_{4}}^{a_{7}}$$ 
and 
$$\coprod_{a_{1}, a_{2}, a_{3}, a_{4}, a_{5}, a_{7}
\in {\cal A}}\left(\coprod_{a_{6}\in {\cal A}}{\cal
V}_{a_{1}a_{2}}^{a_{6}}\otimes {\cal V}_{a_{6}a_{7}}^{a_{5}}\right)
\otimes {\cal V}_{a_{3}a_{4}}^{a_{7}},$$ 
respectively.
Thus the fusing isomorphisms 
$$F(a_{1}, a_{2}, a_{7}; a_{5}):
\coprod_{a_{6}\in {\cal A}}{\cal
V}_{a_{1}a_{6}}^{a_{5}}\otimes {\cal V}_{a_{2}a_{7}}^{a_{6}}
\to \coprod_{a_{6}\in {\cal A}}{\cal
V}_{a_{1}a_{2}}^{a_{6}}\otimes {\cal V}_{a_{6}a_{7}}^{a_{5}}$$
for $a_{1}, a_{2}, a_{7}, a_{5}\in {\cal A}$ give an isomorphism 
\begin{eqnarray*}
\lefteqn{F^{(1)}_{12}: \coprod_{a_{1}, a_{2}, a_{3}, a_{4}, 
a_{5}, a_{6}, a_{7}
\in {\cal A}}{\cal
V}_{a_{1}a_{6}}^{a_{5}}\otimes {\cal V}_{a_{2}a_{7}}^{a_{6}}
\otimes {\cal V}_{a_{3}a_{4}}^{a_{7}}}\nno\\
&&\quad\quad\quad\quad\quad\quad\quad\to \coprod_{a_{1}, a_{2}, 
a_{3}, a_{4}, a_{5}, a_{6}, a_{7}
\in {\cal A}}{\cal
V}_{a_{1}a_{2}}^{a_{6}}\otimes {\cal V}_{a_{6}a_{7}}^{a_{5}}
\otimes {\cal V}_{a_{3}a_{4}}^{a_{7}}.
\end{eqnarray*}
Similarly we have isomorphisms 
\begin{eqnarray*}
\lefteqn{F^{(2)}_{12}: \coprod_{a_{1}, a_{2}, a_{3}, 
a_{4}, a_{5}, a_{6}, a_{7}
\in {\cal A}}{\cal
V}_{a_{1}a_{7}}^{a_{6}}\otimes {\cal V}_{a_{2}a_{3}}^{a_{7}}
\otimes {\cal V}_{a_{}a_{4}}^{a_{5}}}\nno\\
&&\quad\quad\quad\quad\quad\quad\quad
\to \coprod_{a_{1}, a_{2}, a_{3}, a_{4}, a_{5}, a_{6}, a_{7}
\in {\cal A}}{\cal
V}_{a_{1}a_{2}}^{a_{6}}\otimes {\cal V}_{a_{6}a_{3}}^{a_{7}}
\otimes {\cal V}_{a_{7}a_{4}}^{a_{5}},
\end{eqnarray*}
\begin{eqnarray*}
\lefteqn{F_{13}: \coprod_{a_{1}, a_{2}, a_{3}, a_{4}, a_{5}, a_{6}, a_{7}
\in {\cal A}}{\cal
V}_{a_{1}a_{6}}^{a_{5}}\otimes {\cal
V}_{a_{2}a_{3}}^{a_{7}}\otimes {\cal V}_{a_{7}a_{4}}^{a_{6}}}\nno\\
&&\quad\quad\quad\quad\quad\quad\quad
\to \coprod_{a_{1}, a_{2}, a_{3}, a_{4}, a_{5}, a_{6}, a_{7}
\in {\cal A}}{\cal
V}_{a_{1}a_{7}}^{a_{6}}\otimes {\cal V}_{a_{2}a_{3}}^{a_{7}}
\otimes {\cal V}_{a_{}a_{4}}^{a_{5}},
\end{eqnarray*}
\begin{eqnarray*}
\lefteqn{F^{(1)}_{23}: \coprod_{a_{1}, a_{2}, 
a_{3}, a_{4}, a_{5}, a_{6}, a_{7}
\in {\cal A}}{\cal
V}_{a_{1}a_{6}}^{a_{5}}\otimes {\cal V}_{a_{2}a_{7}}^{a_{6}}
\otimes {\cal V}_{a_{3}a_{4}}^{a_{7}}}\nno\\
&&\quad\quad\quad\quad\quad\quad\quad
\to \coprod_{a_{1}, a_{2}, a_{3}, a_{4}, a_{5}, a_{6}, a_{7}
\in {\cal A}}{\cal
V}_{a_{1}a_{6}}^{a_{5}}\otimes {\cal
V}_{a_{2}a_{3}}^{a_{7}}\otimes {\cal V}_{a_{7}a_{4}}^{a_{6}}
\end{eqnarray*}
and 
\begin{eqnarray*}
\lefteqn{F^{(2)}_{23}: \coprod_{a_{1}, a_{2}, a_{3}, 
a_{4}, a_{5}, a_{6}, a_{7}
\in {\cal A}}{\cal
V}_{a_{1}a_{2}}^{a_{6}}\otimes {\cal V}_{a_{6}a_{7}}^{a_{5}}
\otimes {\cal V}_{a_{3}a_{4}}^{a_{7}}}\nno\\
&&\quad\quad\quad\quad\quad\quad\quad
\to \coprod_{a_{1}, a_{2}, a_{3}, a_{4}, a_{5}, a_{6}, a_{7}
\in {\cal A}}{\cal
V}_{a_{1}a_{2}}^{a_{6}}\otimes {\cal V}_{a_{6}a_{3}}^{a_{7}}
\otimes {\cal V}_{a_{7}a_{4}}^{a_{5}}.
\end{eqnarray*}

Also, {}from $\Omega$ and its inverse, 
we  obtain the following maps in the obvious way:
$$\Omega^{(1)}, (\Omega^{-1})^{(1)}: 
\coprod_{a_{1}, a_{2}, a_{3}, a_{4}, a_{5}\in {\cal A}}
{\cal V}_{a_{1}a_{5}}^{a_{4}}\otimes {\cal V}_{a_{2}a_{3}}^{a_{5}}
\to \coprod_{a_{1}, a_{2}, a_{3}, a_{4}, a_{5}\in {\cal A}}
{\cal V}_{a_{2}a_{3}}^{a_{5}}\otimes {\cal V}_{a_{5}a_{1}}^{a_{4}},$$
$$\Omega^{(2)}, (\Omega^{-1})^{(2)}: \coprod_{a_{1}, a_{2}, 
a_{3}, a_{4}, a_{5}\in {\cal A}}
{\cal V}_{a_{1}a_{2}}^{a_{5}}\otimes {\cal V}_{a_{5}a_{3}}^{a_{4}}
\to \coprod_{a_{1}, a_{2}, a_{3}, a_{4}, a_{5}\in {\cal A}}
{\cal V}_{a_{2}a_{1}}^{a_{5}}\otimes {\cal V}_{a_{5}a_{3}}^{a_{4}},$$
$$\Omega^{(3)}, (\Omega^{-1})^{(3)}: 
\coprod_{a_{1}, a_{2}, 
a_{3}, a_{4}, a_{5}\in {\cal A}}
{\cal V}_{a_{1}a_{2}}^{a_{5}}\otimes {\cal V}_{a_{5}a_{3}}^{a_{4}}
\to \coprod_{a_{1}, a_{2}, a_{3}, a_{4}, a_{5}\in {\cal A}}
{\cal V}_{a_{3}a_{5}}^{a_{4}}\otimes {\cal V}_{a_{1}a_{2}}^{a_{5}},$$
$$\Omega^{(4)}, (\Omega^{-1})^{(4)}: 
\coprod_{a_{1}, a_{2}, a_{3}, a_{4}, 
a_{5}\in {\cal A}}
{\cal V}_{a_{1}a_{5}}^{a_{4}}\otimes {\cal V}_{a_{2}a_{3}}^{a_{5}}
\to \coprod_{a_{1}, a_{2}, a_{3}, a_{4}, a_{5}\in {\cal A}}
{\cal V}_{a_{1}a_{5}}^{a_{4}}\otimes {\cal V}_{a_{3}a_{2}}^{a_{5}}.$$

Applying associativity, commutativity and the skew-symmetry of
intertwining operators to products or iterates of three intertwining
operators, we can verify easily that the skew-symmetry and fusing
isomorphisms satisfy the following {\it genus-zero Moore-Seiberg
equations} which were first derived (in the sense explained below) by
Moore and Seiberg \cite{MS}:
\begin{eqnarray}
F^{(2)}_{23}\circ F^{(1)}_{12}&=&F^{(2)}_{12}\circ F_{13}
\circ F_{23}^{(1)},\label{pentagon}\\
F\circ \Omega^{(3)} \circ F&=&\Omega^{(2)}\circ
F\circ \Omega^{(4)},\label{hexagon1}\\
F\circ (\Omega^{-1})^{(3)} \circ F&=&(\Omega^{-1})^{(2)}\circ
F\circ (\Omega^{-1})^{(4)}.\label{hexagon2}
\end{eqnarray}
It is important to note that when Moore and Seiberg derived these equations,
they explicitly assumed the duality properties 
(the convergence property,
associativity and commutativity) of chiral vertex operators (intertwining 
operators). These properties are actually much stronger than the 
equations (\ref{pentagon})--(\ref{hexagon2}) and were  proved 
for intertwining operators for 
suitable vertex operator algebras in \cite{H2}, \cite{H3} and \cite{HL7}.

We call (\ref{pentagon}) the {\it pentagon identity} 
and (\ref{hexagon1}) and 
(\ref{hexagon2}) the {\it hexagon identities} since they correspond to 
the commutativity of the pentagon and hexagon diagrams
for braided tensor categories. Since the formulation of the Jacobi
identity and the proof that intertwining operators for a suitable
vertex operator algebra satisfy the Jacobi identity
do not use these identities, we shall postpone
their proof until the forthcoming paper \cite{H6}.

We now define two maps, corresponding to the multiplication and
iteration of
intertwining operators. The first one is 
\begin{eqnarray*}
{\bf P}: \coprod_{a_{1}, a_{2}, a_{3}, a_{4},
a_{5}\in {\cal A}}{\cal V}_{a_{1}a_{5}}^{a_{4}}\otimes 
{\cal V}_{a_{2}a_{3}}^{a_{5}}&\to &
(\hom(W\otimes W\otimes W, W))\{x_{1}, x_{2}\}\\
{\cal Z}&\mapsto& {\bf P}({\cal Z})
\end{eqnarray*}
defined using products of intertwining operators 
as follows: For 
$${\cal Z}\in  
\coprod_{a_{1}, a_{2}, a_{3}, a_{4},
a_{5}\in {\cal A}}{\cal V}_{a_{1}a_{5}}^{a_{4}}\otimes 
{\cal V}_{a_{2}a_{3}}^{a_{5}},$$ 
the element ${\bf P}({\cal Z})$ to be defined 
can also be
viewed as a map {}from $W\otimes W \otimes W$ to 
$W\{x_{1}, x_{2}\}$. For any $w_{1}, w_{2}, w_{3}\in W$,
we denote the image of $w_{1}\otimes
w_{2}\otimes w_{3}$ under this map by 
$({\bf P}({\cal Z}))(w_{1},  w_{2}, w_{3}; x_{1}, x_{2})$.
Then we define
\begin{eqnarray*}
\lefteqn{({\bf P}({\cal Y}_{1}\otimes 
{\cal Y}_{2}))(w_{(a_{6})},  w_{(a_{7})}, w_{(a_{8})}; x_{1}, x_{2})}\nno\\
&&=\left\{\begin{array}{ll}
{\cal Y}_{1}(w_{(a_{6})}, x_{1}){\cal Y}_{2}(w_{(a_{7})}, x_{2})w_{(a_{8})}&
a_{6}=a_{1}, a_{7}=a_{2}, a_{8}=a_{3}\\
0&\mbox{\rm otherwise}
\end{array}\right.
\end{eqnarray*}
for $a_{1}, \dots, a_{8}\in {\cal A}$, ${\cal Y}_{1}\in {\cal
V}_{a_{1}a_{5}}^{a_{4}}$, ${\cal Y}_{2}\in {\cal
V}_{a_{2}a_{3}}^{a_{5}}$, and $w_{(a_{7})}\in W_{a_{6}}$,
$w_{(a_{7})}\in W_{a_{7}}$, $w_{(a_{8})}\in W_{a_{8}}$.

The second map is
\begin{eqnarray*}
{\bf I}: \coprod_{a_{1}, a_{2}, a_{3}, a_{4},
a_{5}\in {\cal A}}{\cal V}_{a_{1}a_{2}}^{a_{5}}\otimes {\cal
V}_{a_{5}a_{3}}^{a_{4}}&\to&
(\hom(W\otimes W\otimes W, W))\{x_{0}, x_{2}\}\\
{\cal Z}&\mapsto& {\bf I}({\cal Z})
\end{eqnarray*}
defined similarly using iterates of intertwining operators as follows: 
For 
$${\cal Z}\in \coprod_{a_{1}, a_{2}, a_{3}, a_{4},
a_{5}\in {\cal A}}{\cal V}_{a_{1}a_{2}}^{a_{5}}\otimes {\cal
V}_{a_{5}a_{3}}^{a_{4}},$$
the element ${\bf I}({\cal Z})$ to be defined 
can be
viewed as a map {}from $W\otimes W \otimes W$ to 
$W\{x_{0}, x_{2}\}$. For any $w_{1}, w_{2}, w_{3}\in W$,
we denote the image of $w_{1}\otimes
w_{2}\otimes w_{3}$ under this map by 
$({\bf I}({\cal Z}))(w_{(a_{6})}, w_{(a_{7})}, w_{(a_{6})}; 
x_{0}, x_{2})$. We define
\begin{eqnarray*}
\lefteqn{({\bf I}({\cal Y}_{1}\otimes 
{\cal Y}_{2}))(w_{(a_{6})}, w_{(a_{7})}, w_{(a_{8})}; x_{0}, x_{2})}\nno\\
&&=\left\{\begin{array}{ll}
{\cal Y}_{2}({\cal Y}_{1}(w_{(6)}, x_{0})w_{(7)}, x_{2})w_{(a_{8})}
&a_{6}=a_{1}, a_{7}=a_{2}, a_{8}=a_{3}\\
0&\mbox{\rm otherwise}
\end{array}\right.
\end{eqnarray*}
for $a_{1}, \dots, a_{8}\in {\cal A}$, ${\cal Y}_{1}\in {\cal
V}_{a_{1}a_{2}}^{a_{5}}$, ${\cal Y}_{2}\in {\cal
V}_{a_{5}a_{3}}^{a_{4}}$, and $w_{(a_{6})}\in W_{a_{6}}$,
$w_{(a_{7})}\in W_{a_{7}}$, $w_{(a_{8})}\in W_{a_{8}}$. 

We shall call 
${\bf P}$ and ${\bf I}$ the {\it multiplication 
of intertwining operators} and
the {\it iteration of intertwining operators}, respectively. 

We also need to discuss certain special multivalued analytic
functions. Consider the simply connected regions in ${\Bbb C}^{2}$
obtained by cutting the regions $|z_{1}|>|z_{2}|>0$,
$|z_{2}|>|z_{1}|>0$ and $|z_{2}|>|z_{1}-z_{2}|>0$ along the
intersections of these regions with ${\Bbb R}^{2}\subset {\Bbb
C}^{2}$. We denote them by $R_{1}$, $R_{2}$ and $R_{3}$, respectively.
For $a_{1}, a_{2}, a_{3}, a_{4}\in {\cal A}$, let ${\Bbb G}^{a_{1},
a_{2}, a_{3}, a_{4}}$ be the space of multivalued analytic functions on
$$M^{2}=\{(z_{1}, z_{2})\in {\Bbb C}^{2}\;|\; z_{1}, z_{2}\ne 0,
z_{1}\ne z_{2}\}$$
with a choice of a single
valued branch on  the region $R_{1}$ satisfying the following property:
In the regions $|z_{1}|>|z_{2}|>0$,
$|z_{2}|>|z_{1}|>0$ and $|z_{2}|>|z_{1}-z_{2}|>0$, any (multivalued)
branch of $f(z_{1}, z_{2})\in {\Bbb G}^{a_{1}, a_{2}, a_{3}, a_{4}}$ 
can be expanded as
\begin{eqnarray*}
&{\displaystyle \sum_{a\in {\cal A}}z_{1}^{h_{a_{4}}-h_{a_{1}}-h_{a}}
z_{2}^{h_{a}-h_{a_{2}}-h_{a_{3}}}F_{a}(z_{1}, z_{2})},&\\
&{\displaystyle \sum_{a\in {\cal A}}z_{2}^{h_{a_{4}}-h_{a_{2}}-h_{a}}
z_{1}^{h_{a}-h_{a_{1}}-h_{a_{3}}}G_{a}(z_{1}, z_{2})}&
\end{eqnarray*}
and 
$$\sum_{a\in {\cal A}}z_{2}^{h_{a_{4}}-h_{a}-h_{a_{3}}}
(z_{1}-z_{2})^{h_{a}-h_{a_{1}}-h_{a_{2}}}H_{a}(z_{1}, z_{2}),$$
respectively, where for $a\in {\cal A}$, 
$$F_{a}(z_{1}, z_{2})\in {\Bbb C}[z_{1}, z_{1}^{-1}, z_{2}, z_{2}^{-1}]
[[z_{2}/z_{1}]],$$
$$H_{a}(z_{2}, z_{1})\in {\Bbb C}[z_{1}, z_{1}^{-1}, z_{2}, z_{2}^{-1}]
[[z_{1}/z_{2}]]$$
and 
$$G_{a}(z_{1}, z_{1}-z_{2})\in {\Bbb C}[z_{1}, z_{1}^{-1}, z_{1}-z_{2}, 
(z_{1}-z_{2})^{-1}]
[[(z_{1}-z_{2})/z_{2}]].$$

We call the chosen single valued branch on $R_{1}$ of an element 
of ${\Bbb G}^{a_{1}, a_{2}, a_{3}, a_{4}}$
the {\it prefered branch on $R_{1}$}.
The values of the prefered branch of an element of 
${\Bbb G}^{a_{1}, a_{2}, a_{3}, a_{4}}$ 
on the set $R_{1}\cap R_{3}\cap {\Bbb R}^{2}$ give a single valued
branch of the element on $R_{3}$. We call this branch the
{\it prefered branch on $R_{3}$}. Similarly, the values of the
prefered branch on the set $R_{3}\cap R_{2}\cap {\Bbb R}^{2}$ give a
single valued branch of the element on $R_{2}$ and we
call this branch the {\it prefered branch on $R_{2}$}.

Given an element of 
${\Bbb G}^{a_{1}, a_{2}, a_{3}, a_{4}}$, 
the prefered branches of this function on $R_{1}$, $R_{2}$ and $R_{3}$
give  formal series in 
$$\coprod_{a\in {\cal A}}x_{1}^{h_{a_{4}}-h_{a_{1}}-h_{a}}
x_{2}^{h_{a}-h_{a_{2}}-h_{a_{3}}}
{\Bbb C}[x_{1}, x_{1}^{-1}, x_{2}, x_{2}^{-1}][[x_{2}/x_{1}]],$$
$$\coprod_{a\in {\cal A}}x_{2}^{h_{a_{4}}-h_{a_{2}}-h_{a}}
x_{1}^{h_{a}-h_{a_{1}}-h_{a_{3}}}
{\Bbb C}[x_{1}, x_{1}^{-1}, x_{2}, x_{2}^{-1}][[x_{2}/x_{1}]]$$
and 
$$\coprod_{a\in {\cal A}}x_{2}^{h_{a_{4}}-h_{a}-h_{a_{3}}}
x_{0}^{h_{a}-h_{a_{1}}-h_{a_{2}}}{\Bbb C}[
x_{0}, x_{0}^{-1}, x_{2}, x_{2}^{-1}][[x_{0}/x_{2}]],$$
respectively. Thus we have linear maps 
\begin{eqnarray*}
\iota_{12}: {\Bbb G}^{a_{1}, a_{2}, a_{3}, a_{4}}&\to& 
\coprod_{a\in {\cal A}}x_{1}^{h_{a_{4}}-h_{a_{1}}-h_{a}}
x_{2}^{h_{a}-h_{a_{2}}-h_{a_{3}}}{\Bbb C}[x_{1}, x_{1}^{-1}, x_{2},
x_{2}^{-1}][[x_{2}/x_{1}]]\\
\iota_{21}: {\Bbb G}^{a_{1}, a_{2}, a_{3}, a_{4}}&\to& 
\coprod_{a\in {\cal A}}x_{2}^{h_{a_{4}}-h_{a_{2}}-h_{a}}
x_{1}^{h_{a}-h_{a_{1}}-h_{a_{3}}}{\Bbb C}[x_{1}, x_{1}^{-1}, x_{2}, 
x_{2}^{-1}][[x_{2}/x_{1}]]\\
\iota_{20}: {\Bbb G}^{a_{1}, a_{2}, a_{3}, a_{4}}&\to& 
\coprod_{a\in {\cal A}}x_{2}^{h_{a_{4}}-h_{a}-h_{a_{3}}}
x_{0}^{h_{a}-h_{a_{1}}-h_{a_{2}}}{\Bbb C}[ x_{0}, 
x_{0}^{-1}, x_{2}, x_{2}^{-1}][[x_{0}/x_{2}]]
\end{eqnarray*}
generalizing $\iota_{12}$, $\iota_{21}$ and $\iota_{20}$ discussed
before. Since analytic extensions are unique, these maps are
injective.

For any $a_{1}, a_{2}, a_{3}, a_{4}\in {\cal A}$, 
${\Bbb G}^{a_{1}, a_{2}, a_{3}, a_{4}}$ is a module over the ring
$${\Bbb C}[x_{1}, x_{1}^{-1}, x_{2}, x_{2}^{-1}, (x_{1}-x_{2})^{-1}].$$
We have:

\begin{lemma}\label{free}
The module ${\Bbb G}^{a_{1}, a_{2}, a_{3}, a_{4}}$ is free.
\end{lemma}

The proof of this lemma is given in the next section.
We fix a basis $\{f^{a_{1}, a_{2}, a_{3}, a_{4}}_{\alpha}\}_{\alpha\in
{\Bbb A}}$ of this free module ${\Bbb G}^{a_{1}, a_{2}, a_{3}, a_{4}}$.

\begin{theo}\label{jacobi}
For any $a_{1}, a_{2}, a_{3}, a_{4}\in {\cal A}$, there exist 
linear maps
\begin{eqnarray}\label{1.-4}
\lefteqn{g^{a_{1}, a_{2}, a_{3},
a_{4}}_{\alpha}: W_{a_{1}}\otimes W_{a_{2}}\otimes W_{a_{3}}\otimes 
\coprod_{a_{5}\in {\cal A}}
{\cal V}_{a_{1}a_{5}}^{a_{4}}\otimes 
{\cal V}_{a_{2}a_{3}}^{a_{5}}}\nno\\
&&\hspace{8em}\to  W_{a_{4}}[x_{1}, x_{1}^{-1}, x_{2},
x_{2}^{-1}][[x_{2}/x_{1}]]\nno\\
&&\hspace{3em} w_{(a_{1})}\otimes w_{(a_{2})}\otimes w_{(a_{3})}\otimes {\cal
Z}\nno\\
&&\hspace{8em} \mapsto g^{a_{1}, a_{2}, a_{3},
a_{4}}_{\alpha}(w_{(a_{1})},
w_{(a_{2})}, w_{(a_{3})}, {\cal Z}; x_{1}, x_{2})
\end{eqnarray}
and 
\begin{eqnarray}\label{1.-3}
\lefteqn{h^{a_{1}, a_{2}, a_{3},
a_{4}}_{\alpha}: W_{a_{1}}\otimes W_{a_{2}}\otimes W_{a_{3}}\otimes 
\coprod_{a_{5}\in {\cal A}}
{\cal V}_{a_{1}a_{2}}^{a_{5}}\otimes 
{\cal V}_{a_{5}a_{3}}^{a_{4}}}\nno\\
&&\hspace{8em}\to W_{a_{4}}[x_{0}, 
x_{0}^{-1}, x_{2}, x_{2}^{-1}][[x_{0}/x_{2}]]\nno\\
&&\hspace{3em} w_{(a_{1})}\otimes w_{(a_{2})}\otimes w_{(a_{3})}\otimes {\cal
Z}\nno\\
&&\hspace{8em} \mapsto h^{a_{1}, a_{2}, a_{3},
a_{4}}_{\alpha}(w_{(a_{1})},
w_{(a_{2})}, w_{(a_{3})}, {\cal Z}; x_{0}, x_{2})
\end{eqnarray}
for $\alpha\in {\Bbb A}$,
such that  for
any
$w_{(a_{1})}\in W_{a_{1}}$, $w_{(a_{2})}\in W_{a_{2}}$, $w_{(a_{3})}\in 
W_{a_{3}}$, and any
$${\cal Z}\in \coprod_{a_{4}, a_{5}\in {\cal A}}
{\cal V}_{a_{1}a_{5}}^{a_{4}}\otimes 
{\cal V}_{a_{2}a_{3}}^{a_{5}}\subset \coprod_{a_{1}, a_{2}, a_{3}, a_{4},
a_{5}\in {\cal A}}{\cal V}_{a_{1}a_{5}}^{a_{4}}\otimes 
{\cal V}_{a_{2}a_{3}}^{a_{5}},$$
only finitely many of $$g^{a_{1}, a_{2}, a_{3},
a_{4}}_{\alpha}(w_{(a_{1})},
w_{(a_{2})}, w_{(a_{3})}, {\cal Z}; x_{1}, x_{2})$$ and 
$$h^{a_{1}, a_{2}, a_{3},
a_{4}}_{\alpha}(w_{(a_{1})},
w_{(a_{2})}, w_{(a_{3})}, {\cal F}({\cal Z}); x_{0}, x_{2}),$$
$\alpha\in {\Bbb A}$, are nonzero,
\begin{eqnarray}\label{1.-2}
\lefteqn{({\bf P}({\cal Z}))(w_{(a_{1})},
w_{(a_{2})}, w_{(a_{3})}; x_{1}, x_{2})}\nno\\
&&=\sum_{\alpha\in {\Bbb A}}g^{a_{1}, a_{2}, a_{3},
a_{4}}_{\alpha}(w_{(a_{1})},
w_{(a_{2})}, w_{(a_{3})}, {\cal Z}; x_{1}, x_{2})\iota_{12}\left(f^{a_{1}, 
a_{2}, a_{3}, a_{4}}_{\alpha}\right),
\end{eqnarray}
\begin{eqnarray}\label{1.-1}
\lefteqn{
({\bf I}({\cal F}({\cal Z})))(w_{(a_{1})}, w_{(a_{2})}, 
w_{(a_{3})}; x_{0}, x_{2})}\nno\\
&&=\sum_{\alpha\in {\Bbb A}}h^{a_{1}, a_{2}, a_{3},
a_{4}}_{\alpha}(w_{(a_{1})},
w_{(a_{2})}, w_{(a_{3})}, {\cal F}({\cal Z}); x_{0}, x_{2})\iota_{20}
(f^{a_{1}, 
a_{2}, a_{3}, a_{4}}_{\alpha}),
\end{eqnarray}
and the following {\it Jacobi identity} holds:
\begin{eqnarray}\label{1.0}
\lefteqn{x_{0}^{-1}
\delta\left(\frac{x_{1}-x_{2}}{x_{0}}\right)
g^{a_{1}, a_{2}, a_{3},
a_{4}}_{\alpha}(w_{(a_{1})},
w_{(a_{2})}, w_{(a_{3})}, {\cal Z}; x_{1}, x_{2})}\nno\\
&&\quad -x_{0}^{-1}\delta\left(\frac{x_{2}-x_{1}}{-x_{0}}\right)
g^{a_{2}, a_{1}, a_{3},
a_{4}}_{\alpha}(w_{(a_{2})},
w_{(a_{1})}, w_{(a_{3})}, {\cal B}({\cal Z}); x_{2}, x_{1})
\nno\\
&&=x_{2}^{-1}\delta\left(\frac{x_{1}-x_{0}}{x_{2}}\right)
h^{a_{1}, a_{2}, a_{3},
a_{4}}_{\alpha}(w_{(a_{1})},
w_{(a_{2})}, w_{(a_{3})}, {\cal F}({\cal Z}); x_{0}, x_{2})
\end{eqnarray}
for $\alpha\in {\Bbb A}$.
\end{theo}

Theorem \ref{jacobi} will be proved in Section 4 based on the
generalized
rationality and duality property in formal variables proved in the
same section. 
All the properties of the Jacobi identity for vertex operator algebras
have generalizations for this Jacobi identity. This Jacobi identity
can be used to replace the associativity and  skew-symmetry in the
definition of intertwining operator algebra. 
These will be discussed in the forthcoming paper \cite{H6} on the axiomatic
properties of intertwining operator algebras.

The Jacobi identity (\ref{1.0}) is simple but sometimes 
is  not explicit enough for
concrete calculations. To obtain an explicit form,
we choose a basis ${\cal
Y}^{a_{3}}_{a_{1}a_{2};i}$ of the space ${\cal
V}_{a_{1}a_{2}}^{a_{3}}$  for $a_{1}, a_{2}, a_{3}\in {\cal A}$ 
and $i=1, \dots, {\cal
N}_{a_{1}a_{2}}^{a_{3}}$. Then the 
{\it fusing matrices} under this basis 
corresponding to the fusing isomorphisms 
${\cal F}(a_{1},
a_{2}, a_{3}; a_{4})$, $a_{1}, a_{2}, a_{3}, a_{4}\in {\cal A}$,
whose entries
are given by
\begin{eqnarray*}
\lefteqn{({\cal F}(a_{1},
a_{2}, a_{3}; a_{4}))({\cal Y}^{a_{4}}_{a_{1}a_{5};i}\otimes
{\cal Y}^{a_{5}}_{a_{2}a_{3};j})}\nno\\
&&=\sum_{a\in {\cal
A}}\sum_{k=1}^{{\cal N}_{a_{1}a_{2}}^{a}}
\sum_{l=1}^{{\cal N}_{aa_{3}}^{a_{4}}}
{\cal F}_{a_{5};  a}^{i, j; k, l}(a_{1}, a_{2}, a_{3}; a_{4})({\cal
Y}_{a_{1}a_{2}; k}^{a}\otimes {\cal Y}_{aa_{3}; l}^{a_{4}})
\end{eqnarray*}
for $a_{5}\in {\cal A}$, $i=1, \dots, {\cal
N}_{a_{1}a_{5}}^{a_{4}}$, $j=1, \dots, {\cal
N}_{a_{2}a_{3}}^{a_{5}}$.

Using the same basis, we obtain the
{\it braiding matrices} corresponding to the braiding isomorphisms
${\cal B}(a_{1},
a_{2}; a_{3}, a_{4})$, $a_{1}, a_{2}, a_{3}\in {\cal A}$,
whose entries are given by
\begin{eqnarray*}
\lefteqn{({\cal B}(a_{1},
a_{2}; a_{3}, a_{4}))({\cal Y}^{a_{4}}_{a_{1}a_{5};i}\otimes
{\cal Y}^{a_{5}}_{a_{2}a_{3}; j})}\nno\\
&&=\sum_{a\in {\cal A}}\sum_{k=1}^{{\cal N}_{a_{2}a}^{a_{4}}}
\sum_{l=1}^{{\cal N}_{a_{1}a_{3}}^{a}}
{\cal B}_{a_{5}; a}^{i, j; k, l}(a_{1}, a_{2}; a_{3}, a_{4})({\cal
Y}_{a_{2}a; k}^{a_{4}}\otimes {\cal Y}_{a_{1}a_{3}; l}^{a})
\end{eqnarray*}
for $a_{5}\in {\cal A}$, $i=1, \dots, {\cal
N}_{a_{1}a_{5}}^{a_{4}}$, $j=1, \dots, {\cal N}_{a_{2}a_{3}}^{a_{5}}$.

Using the braiding and fusing matrices above, the Jacobi identity
can be written as
\begin{eqnarray}\label{1.1}
\lefteqn{x_{0}^{-1}
\delta\left(\frac{x_{1}-x_{2}}{x_{0}}\right)
g^{a_{1}, a_{2}, a_{3}, a_{4}}_{\alpha}(w_{(a_{1})},
w_{(a_{2})}, w_{(a_{3})}, {\cal
Y}_{a_{1}a_{5}; i}^{a_{4}}\otimes {\cal Y}_{a_{2}a_{3}; j}^{a_{5}}; 
x_{1}, x_{2})}\nno\\
&&\quad -
x_{0}^{-1}\delta\left(\frac{x_{2}-x_{1}}{-x_{0}}\right)
\sum_{a\in {\cal
A}}\sum_{k=1}^{{\cal N}_{a_{2}a}^{a_{4}}} \sum_{l=1}^{{\cal
N}_{a_{1}a_{3}}^{a}} {\cal B}_{a_{5}; a}^{i, j; k, l}(a_{1}, 
a_{2}; a_{3}, a_{4})\cdot\nno\\
&&\hspace{5em}\cdot
g^{a_{2}, a_{1}, a_{3}, a_{4}}
_{\alpha}(w_{(a_{2})},
w_{(a_{1})}, w_{(a_{3})}, {\cal Y}_{a_{2}a; k}^{a_{4}}
\otimes {\cal
Y}_{a_{1}a_{3}; l}^{a}; 
x_{2}, x_{1}))\nno\\
&&=x_{2}^{-1}\delta\left(\frac{x_{1}-x_{0}}{x_{2}}\right)
 \sum_{a\in {\cal
A}}\sum_{k=1}^{{\cal N}_{a_{1}a_{2}}^{a}} \sum_{l=1}^{{\cal
N}_{aa_{3}}^{a_{4}}} {\cal F}_{a_{5}; a}^{i, j; k, l}(a_{1}, a_{2}, a_{3};
a_{4})\cdot\nno\\
&&\hspace{5em}\cdot 
h^{a_{1}, a_{2}, a_{3}, a_{4}}
_{\alpha}(w_{(a_{1})},
w_{(a_{2})}, w_{(a_{3})}, {\cal Y}_{a_{2}a; k}^{a_{4}}
\otimes {\cal
Y}_{a_{1}a_{3}; l}^{a}; 
x_{0}, x_{2}),
\end{eqnarray}
for $a_{1}, a_{2}, a_{3}, a_{4}, a_{5}\in {\cal A}$,  
 $i=1, \dots,  {\cal N}_{a_{1}a_{5}}^{a_{4}}$, $j=1, \dots, {\cal
N}_{a_{2}a_{3}}^{a_{5}}$, 
$w_{(a_{1})}\in W_{a_{1}}$, $w_{(a_{2})}\in W_{a_{2}}$, and 
$\alpha\in {\Bbb A}$.

\renewcommand{\theequation}{\thesection.\arabic{equation}}
\renewcommand{\therema}{\thesection.\arabic{rema}}
\setcounter{equation}{0}
\setcounter{rema}{0}

\section{Generalized rationality, duality and the 
proof of Theorem \ref{jacobi}}

We first reformulate the associativity and commutativity 
for intertwining operators using
the basis ${\cal
Y}^{a_{3}}_{a_{1}a_{2};i}$,  $a_{1}, a_{2}, a_{3}\in {\cal A}$,
$i=1, \dots, {\cal
N}_{a_{1}a_{2}}^{a_{3}}$, chosen in the preceding section
and the fusing and braiding matrices as follows:

\begin{description}

\item[Associativity:] For any $a_{1}, a_{2}, a_{3}, a_{4},
a_{5}\in {\cal A}$, any $i, j\in {\Bbb Z}$ satisfying 
$1\le i\le {\cal N}_{a_{1}a_{5}}^{a_{4}}$, 
$1\le j\le {\cal N}_{a_{2}a_{3}}^{a_{5}}$, 
the (multivalued) analytic function
\begin{equation}\label{2.1}
\langle w'_{(a_{4})}, 
{\cal Y}_{a_{1}a_{5}; i}^{a_{4}}(w_{(a_{1})}, x_{1})
{\cal Y}_{a_{2}a_{3}; j}^{a_{5}}(w_{(a_{2})}, 
x_{2})w_{(a_{3})}\rangle_{W_{a_{4}}}\mbar_{x_{1}=z_{1},
x_{2}=z_{2}}
\end{equation}
defined on the region $|z_{1}|>|z_{2}|>0$ 
and the (multivalued) analytic function
\begin{eqnarray}\label{2.2}
\lefteqn{\sum_{a\in {\cal A}}\sum_{k=1}^{{\cal N}_{a_{1}a_{2}}^{a}}
\sum_{l=1}^{{\cal N}_{aa_{3}}^{a_{4}}}
{\cal F}_{a_{5}; a}^{i, j;k, l}(a_{1}, a_{2}, a_{3}; a_{4})\cdot}\nno\\
&&\cdot\langle w'_{(a_{4})}, {\cal Y}_{aa_{3}; l}^{a_{4}}
({\cal Y}_{a_{1}a_{2}; k}^{a}(w_{1},
x_{0})w_{2}, x_{2})w_{3}\rangle_{W_{a_{4}}}\mbar_{x_{0}=z_{1}-z_{2},
x_{2}=z_{2}}
\end{eqnarray} 
defined on the region
$|z_{2}|>|z_{1}-z_{2}|>0$ are equal on the intersection
$|z_{1}|> |z_{2}|>|z_{1}-z_{2}|>0$. 
In addition,
\begin{eqnarray*}
\lefteqn{\langle w'_{(a_{4})}, 
{\cal Y}_{a_{1}a_{5}; i}^{a_{4}}(w_{(a_{1})}, x_{1})
{\cal Y}_{a_{2}a_{3}; j}^{a_{5}}(w_{(a_{2})}, 
x_{2})w_{(a_{3})}\rangle_{W_{a_{4}}}\mbar_{x_{1}^{n}=e^{n\log z_{1}},
x_{2}^{n}=e^{n\log z_{2}}}}\nno\\
&&=\sum_{a\in {\cal A}}\sum_{k=1}^{{\cal N}_{a_{1}a_{2}}^{a}}
\sum_{l=1}^{{\cal N}_{aa_{3}}^{a_{4}}}
{\cal F}_{a_{5}; a}^{i, j;k, l}(a_{1}, a_{2}, a_{3}; a_{4})\cdot\nno\\
&&\quad \cdot
\langle w'_{(a_{4})}, {\cal Y}_{aa_{3}; l}^{a_{4}}
({\cal Y}_{a_{1}a_{2}; k}^{a}(w_{1},
x_{0})w_{2}, x_{2})w_{3}\rangle_{W_{a_{4}}}
\mbar_{x_{0}^{n}=e^{n\log (z_{1}-z_{2})},
x_{2}^{n}=e^{n\log z_{2}}}
\end{eqnarray*}
when $z_{1}$ and $z_{2}$ are positive real numbers and satisfying 
$z_{1}>z_{2}>z_{1}-z_{2}>0$.

\end{description}

\begin{description}

\item[Commutativity:] For any $a_{1}, a_{2}, a_{3}, a_{4},
a_{5}\in {\cal A}$, any $i, j\in {\Bbb Z}$ satisfying 
$1\le i\le {\cal N}_{a_{1}a_{5}}^{a_{4}}$, 
$1\le j\le {\cal N}_{a_{2}a_{3}}^{a_{5}}$,
the (multivalued) analytic function
(\ref{2.1})
on $\{(z_{1}, z_{2})\in {\Bbb C}\times {\Bbb C}\;|\; |z_{1}|>|z_{2}|>0\}$ 
and the (multivalued) analytic function
\begin{eqnarray}\label{2.3}
\lefteqn{\sum_{a\in {\cal A}}\sum_{k=1}^{{\cal N}_{a_{2}a}^{a_{4}}}
\sum_{l=1}^{{\cal N}_{a_{1}a_{3}}^{a}}
{\cal B}_{a_{5}; a}^{i, j; k, l}(a_{1}, a_{2}; a_{3}, a_{4})\cdot}\nno\\
&&\cdot\langle w'_{(a_{4})}, {\cal Y}_{a_{2}a; k}^{a_{4}}
(w_{2}, x_{2}){\cal Y}_{a_{1}a_{3}; k}^{a}(w_{1},
x_{1})w_{3}\rangle_{W_{a_{4}}}\mbar_{x_{1}=z_{1},
x_{2}=z_{2}}
\end{eqnarray}
on
$\{(z_{1}, z_{2})\in {\Bbb C}\times {\Bbb C}\;|\;
|z_{2}|>|z_{1}|>0\}$ are analytic extensions of each other.

\end{description}

We first prove the following:

\begin{lemma}\label{correl1}
For any $a_{1}, a_{2}, a_{3}, a_{4},
a_{5}\in {\cal A}$ and any $i, j\in {\Bbb Z}$ satisfying 
$1\le i\le {\cal N}_{a_{1}a_{5}}^{a_{4}}$, 
$1\le j\le {\cal N}_{a_{2}a_{3}}^{a_{5}}$, respectively, there exists
a multivalued analytic function defined on $M^{2}=\{(z_{1},
z_{2})\in {\Bbb C}^{2}\;|\;z_{1}, z_{2}\ne 0, z_{1}\ne
z_{2}\}$ such that (\ref{2.1}), (\ref{2.2}) and (\ref{2.3}) 
are restrictions of this function to their domains.
\end{lemma}
\pf
The manifold $M^{2}$ is the union of the following parts:

\begin{enumerate}

\item The region given by $|z_{1}|>|z_{2}|>0$.\label{region1}

\item The region given by $|z_{2}|>|z_{1}-z_{2}|>0$.\label{region2}

\item The region given by $|z_{2}|>|z_{1}|>0$.\label{region3}

\item The region given by $|z_{1}|>|z_{1}-z_{2}|>0$.\label{region4}

\item The region given by
$|z_{1}-z_{2}|>|z_{1}|>0$.\label{region5}

\item The region given by
$|z_{1}-z_{2}|>|z_{2}|>0$.\label{region6}

\item The set given by $|z_{1}|=|z_{2}|=|z_{1}-z_{2}|>0$.\label{region7}

\end{enumerate}

By associativity and commutativity, 
we already have a multivalued analytic function of $z_{1}$ and $z_{2}$
defined on the union of the regions \ref{region1},  \ref{region2} and 
\ref{region3}
such that its restrictions to the regions \ref{region1},
\ref{region2} and  \ref{region3}
are indeed (\ref{2.1}), (\ref{2.2}) and (\ref{2.3}), respectively. 
We now extend the 
function to the regions \ref{region4},  \ref{region5}, 
 \ref{region6} and the set \ref{region7}.

Using associativity, we can extend the function we already have to 
region \ref{region4} {}from the function already defined on 
 region \ref{region3}. 
To extend the function to region \ref{region6}, we note that 
by the convergence property,
the sums of the series
\begin{eqnarray*}
\lefteqn{\langle e^{z_{2}L'(1)}w'_{(a_{4})}, 
{\cal Y}_{a_{1}a_{5}; i}^{a_{4}}(w_{(a_{1})}, x_{1})\cdot}\nno\\
&&\hspace{4em}\cdot
(\Omega(a_{1}, a_{2}; a_{3})({\cal Y}_{a_{2}a_{3}; j}^{a_{5}}))(w_{(a_{3})}, 
x_{2})w_{(a_{2})}\rangle_{W_{a_{4}}}\mbar_{x_{1}=z_{1}-z_{2},
x_{2}=-z_{2}}
\end{eqnarray*}
for all choices of powers of $z_{1}-z_{2}$ and $-z_{2}$ 
define a multivalued 
analytic function on  region \ref{region6}.
Also note that the intersection of 
regions \ref{region6} and 
 \ref{region1} is not empty. 
{}From the definition, the skew-symmetry and the conjugation formula 
\begin{equation}\label{conj}
e^{xL(-1)}{\cal Y}_{a_{1}a_{5}; i}^{a_{4}}
(w_{(a_{1})}, x_{1})e^{-xL(-1)}={\cal Y}(w_{(a_{1})}, x_{1}+x)
\end{equation}
for ${\cal Y}_{a_{1}a_{5}; i}^{a_{4}}(w_{(a_{1})}, x_{1})$ which is an
easy consequence of the $L(-1)$-derivative property, 
we see that on the intersection
of  regions \ref{region6} and  \ref{region1}, 
this function is equal to the function already defined on  
region \ref{region1}. 
Thus we obtain the extension we need. Similarly we can obtain the extension
of the function to region \ref{region5}. 

Finally we have to extend the function to  set \ref{region7}. 
Let $(z_{1}^{(0)}, z_{2}^{(0)})\in M^{2}$ be a point in  set \ref{region7}. 
Then we have 
$|z_{1}^{(0)}|=|z_{2}^{(0)}|=|z_{1}^{(0)}-z_{2}^{(0)}|>0$. 
It is easy to see that there always exists a nonzero complex 
number $\epsilon$ such that 
$|z_{1}^{(0)}+\epsilon|>|z_{2}^{(0)}+\epsilon|>|\epsilon|>0$.
We define a 
multivalued analytic function on the region 
$$D_{\epsilon}=\{(z_{1}, z_{2})\in {\Bbb C}^{2}\;|\;
|z_{1}+\epsilon|>|z_{2}+\epsilon|>|\epsilon|>0\}$$
using  the sums of the series
$$
\langle e^{-\epsilon L'(1)}w'_{(a_{4})}, 
{\cal Y}_{a_{1}a_{5}; i}^{a_{4}}(w_{(a_{1})}, x_{1})
{\cal Y}_{a_{2}a_{3}; j}^{a_{5}}(w_{(a_{2})}, 
x_{2})e^{\epsilon L(-1)}w_{(a_{3})}
\rangle_{W_{a_{4}}}\mbar_{x_{1}=z_{1}+\epsilon,
x_{2}=z_{2}+\epsilon}
$$
for all choices of powers of $z_{1}+\epsilon$ and $z_{2}+\epsilon$.
For $z_{1}$ and $z_{2}$ in $D_{\epsilon}$, the series above are indeed
absolutely convergent by the convergence property because using
skew-symmetry, we have
$$e^{\epsilon L(-1)}w_{(a_{3})}
=\Omega(Y_{W^{a_{3}}})(w_{(a_{3})}, x_{3}){\bf 1}\mbar_{x_{3}=\epsilon}$$
where $Y_{W^{a_{3}}}$ 
is the vertex operator map defining the $W^{e}$-structure 
on $W^{a_{3}}$. Thus the series above define
a multivalued analytic function on the region. Since $(z_{1}^{(0)},
z_{2}^{(0)})\in D_{\epsilon}$, this function is defined at
$(z_{1}^{(0)}, z_{2}^{(0)})$. Note that $D_{\epsilon}$ has nonempty
intersection with  region \ref{region1}. Using the conjugation
formula for intertwining operators by $e^{z_{2}L(-1)}$
(cf. (\ref{conj})), 
we see that at
$(z_{1}, z_{2})$ in the intersection of $D_{\epsilon}$ and region
\ref{region1}, the values of the function we just defined is equal to
the values of the function already defined on region \ref{region1}.  So
we obtain an extension of the function to $D_{\epsilon}$. In
particular, we have extended our function to the point $(z_{1}^{(0)},
z_{2}^{(0)})$.  Since this is an arbitrary point in set
\ref{region7}, our function has been extended to set
\ref{region7}.  
\epfv

We need Lemma \ref{free} to formulate the generalized rationality and
duality properties in terms of formal variables.
 
\vspace{1em}

{\it Proof of Lemma \ref{free}}: By definition, elements of 
${\Bbb G}^{a_{1}, a_{2}, a_{3}, a_{4}}$ can be written as linear combinations
of functions of the form 
$z_{1}^{h_{a_{4}}-h_{a_{1}}-h_{a_{2}}-h_{a_{3}}+k}f(z_{2}/z_{1})$ 
where $k\in {\Bbb
Z}$ and
$f(z)$ is a multivalued analytic 
function on ${\Bbb C}\cup \{\infty\}\setminus \{0, 1\}$ with a choice
of a branch on the region obtained by cutting the region $0<|z|<1$
along the open half line $0<z<1$
such that 
in the region $0<|z|<1$, $1<|z|<\infty$ and $0<|z-1|<1$, any 
(multivalued) branches
of $f(z)$ can be
written as 
$\sum_{a\in {\Bbb A}}z^{h_{a}-h_{a_{2}}-h_{a_{3}}}F_{a}(z)$, 
$\sum_{a\in {\Bbb A}}z^{h_{a_{4}}-h_{a_{2}}-h_{a}}G_{a}(z)$ and 
$\sum_{a\in {\Bbb A}}(z-1)^{h_{a}-h_{a_{1}}-h_{a_{2}}}H_{a}(z)$,
respectively, where 
$F_{a}(z)$, $G_{a}(z)$ and $H_{a}(z)$ are analytic functions in the region 
$|z|<1$, $1<|z|$ (including $\infty$) and $|z-1|<1$. 
Now consider the space ${\Bbb F}$ of all such $f(z)$. This space is a 
module over the ring ${\Bbb C}[x,
x^{-1}, (x-1)^{-1}]$. Using Zorn's lemma, we can find a
minimal set of generators of this module. 

We claim that ${\Bbb F}$ is
in fact a free module with a basis given by this set of generators. 
We first prove that ${\Bbb C}[x,
x^{-1}, (x-1)^{-1}]$ is a principal ideal domain. 
Let $I$ be an ideal of ${\Bbb C}[x,
x^{-1}, (x-1)^{-1}]$. Let $I_{P}=I\cap {\Bbb
C}[x]$. Then for any $f\in {\Bbb
C}[x]$ and $g\in I_{P}$, $fg$ is in $I$ and 
in ${\Bbb
C}[x]$. So $fg$ is in $I_{P}$. Thus $I_{P}$ is an ideal of ${\Bbb
C}[x]$. Since ${\Bbb
C}[x_{2}/x_{1}]$ is a principal ideal domain, there exists $f\in
I_{P}$ such that $I_{P}$ is the ideal generated by $f$. Given any
 $g\in I$, there exist $m, n\in {\Bbb N}$ such that 
$(x)^{m}(x-1)^{n}g\in {\Bbb
C}[x]$. Since $I$ is an ideal,
$x^{m}(x-1)^{n}g$ is also in $I$. 
Thus $x^{m}(x-1)^{n}g\in I_{P}$. 
So 
$$x^{m}(x-1)^{n}g=hf$$
for some $h\in {\Bbb
C}[x]$ and equivalently, 
$$g=x^{-m}(x-1)^{-n}hf.$$
Since $g$ is arbitrary, 
we conclude that $I$ is an ideal generated by $f$. Since $I$ is
arbitrary, we see that ${\Bbb C}[x,
x^{-1}, (x-1)^{-1}]$ is a principal ideal domain. 

We now prove the claim that 
${\Bbb F}$ is free and the minimal set of generators is a basis.
If not, there exist $f_{1}, \dots, f_{m}$ in the set of generators 
and $g_{1}, \dots, g_{m}$ in ${\Bbb C}[x,
x^{-1}, (x-1)^{-1}]$ such that 
\begin{equation}\label{ind}
\sum_{i=1}^{m}g_{i}f_{i}=0.  
\end{equation}
Consider the submodule
generated by $f_{1}, \dots, f_{m}$. By the structure theorem for
finitely generated modules over a principal ideal domain, we know that
this submodule is a direct sum of a free submodule and a torsion
submodule. Since  nonzero elements of ${\Bbb
C}[x, x^{-1}, (x-1)^{-1}]$ applied to nonzero
elements of ${\Bbb F}$ are not zero,
there is no torsion submodule. Thus the submodule generated by $f_{1},
\dots, f_{m}$ is free. This contradicts to (\ref{ind}) and the fact
that our set of generators is minimal. So 
${\Bbb F}$ is a free module and the  minimal set of
generators is a basis.

Now we consider functions of the form
$z_{1}^{h_{a_{4}}-h_{a_{1}}-h_{a_{2}}-h_{a_{3}}}f(z_{2}/z_{1})$
where $f(z)$ is a basis element of the free module ${\Bbb
F}$ over ${\Bbb C}[x, x^{-1}, (x-1)^{-1}]$.  Then these functions
generate the module ${\Bbb G}^{a_{1}, a_{2}, a_{3}, a_{4}}$ over the
ring 
$${\Bbb C}[x_{1}, x_{1}^{-1}, x_{2}, x_{2}^{-1},
(x_{1}-x_{2})^{-1}].$$
We prove that these generators must be linearly
independent. Assume that
\begin{equation}\label{ind2}
\sum_{i=1}^{m}g_{i}z_{1}^{h_{a_{4}}-h_{a_{1}}-h_{a_{2}}-h_{a_{3}}}
f_{i}(z_{2}/z_{2})=0
\end{equation}
where for $i=1, \dots, m$,
$g_{i}\in {\Bbb C}[x_{1}, x_{1}^{-1}, x_{2}, x_{2}^{-1},
(x_{1}-x_{2})^{-1}]$ and $f_{i}(z)$ are 
basis elements of the module over ${\Bbb
C}[x, x^{-1}, (x-1)^{-1}]$ discussed above. 
We want to prove that $g_{i}=0$ for $i=1, \dots, m$. 
We write 
$$g_{i}=\sum_{j=n_{1}}^{n_{2}}l_{ij}(x_{2}/x_{1})x_{1}^{j}$$
for $i=1,
\dots, m$, where $l_{ij}$ for  $i=1,
\dots, m$ and $j=n_{1}, \dots, n_{2}$ are elements of ${\Bbb
C}[x, x^{-1}, (x-1)^{-1}]$. 
The equation (\ref{ind2}) becomes 
\begin{equation}\label{ind3}
\sum_{j=n_{1}}^{n_{2}}\sum_{i=1}^{m}l_{ij}(z_{2}/z_{1})
f_{i}(z_{2}/z_{1})z_{1}^{j}=0.
\end{equation}
Since $l_{ij}(z_{2}/z_{1})
f_{i}(z_{2}/z_{1})$ for  $i=1,
\dots, m$ and $j=n_{1}, \dots, n_{2}$ are functions of  variable 
$z_{2}/z_{1}$ which is independent of the variable $z_{1}$, we obtain
$$\sum_{i=1}^{m}l_{ij}(z_{2}/z_{1})
f_{i}(z_{2}/z_{1})=0$$
by (\ref{ind3}), or equivalently 
$$\sum_{i=1}^{m}l_{ij}(x)
f_{i}(z)=0.$$ 
But $f_{i}$ for  $i=1,
\dots, m$ are linearly independent over 
$${\Bbb
C}[x, x^{-1}, (x-1)^{-1}].$$
So 
$l_{ij}=0$ for $i=1,
\dots, m$ and $j=n_{1}, \dots, n_{2}$. Thus $g_{i}=0$ for $i=1,
\dots, m$. 

Since functions of the form
$z^{h_{a_{4}}-h_{a_{1}}-h_{a_{2}}-h_{a_{3}}}_{1}f(z_{2}/z_{1})$
generate the module ${\Bbb G}^{a_{1}, a_{2}, a_{3}, a_{4}}$ and are
linearly independent, ${\Bbb G}^{a_{1}, a_{2}, a_{3}, a_{4}}$ is
free. \epfv

We choose a basis $\{f^{a_{1}, a_{2}, a_{3}, a_{4}}_{\alpha}\}_{\alpha\in
{\Bbb A}}$ of ${\Bbb G}^{a_{1}, a_{2}, a_{3}, a_{4}}$.

We have the following 
generalized rationality of products and commutativity in formal variables:

\begin{theo}\label{rat}
For any $a_{1}, a_{2}, a_{3}, a_{4}\in {\cal A}$, there exist 
linear maps
\begin{eqnarray*}
\lefteqn{g^{a_{1}, a_{2}, a_{3},
a_{4}}_{\alpha}: W_{a_{1}}\otimes W_{a_{2}}\otimes W_{a_{3}}\otimes 
\coprod_{a_{5}\in {\cal A}}
{\cal V}_{a_{1}a_{5}}^{a_{4}}\otimes 
{\cal V}_{a_{2}a_{3}}^{a_{5}}}\nno\\
&&\hspace{8em} \to  W_{a_{4}}[x_{1}, x_{1}^{-1}, x_{2},
x_{2}^{-1}][[x_{2}/x_{1}]]\nno\\
&&\hspace{3em} w_{(a_{1})}\otimes w_{(a_{2})}\otimes w_{(a_{3})}\otimes {\cal
Z}\nno\\
&&\hspace{8em}\mapsto g^{a_{1}, a_{2}, a_{3},
a_{4}}_{\alpha}(w_{(a_{1})},
w_{(a_{2})}, w_{(a_{3})}, {\cal Z}; x_{1}, x_{2}),
\end{eqnarray*}
$\alpha\in {\Bbb A}$,
satisfying the following property: 

{\bf (a) (generalized rationality of products)} For
any
$w_{(a_{1})}\in W_{a_{1}}$, $w_{(a_{2})}\in W_{a_{2}}$, $w_{(a_{3})}\in 
W_{a_{3}}$, $w'_{(a_{4})}\in W'_{a_{4}}$, and any
$${\cal Z}\in \coprod_{a_{5}\in {\cal A}}
{\cal V}_{a_{1}a_{5}}^{a_{4}}\otimes 
{\cal V}_{a_{2}a_{3}}^{a_{5}},$$ 
only finitely many of 
$$g^{a_{1}, a_{2}, a_{3},
a_{4}}_{\alpha}(w_{(a_{1})},
w_{(a_{2})}, w_{(a_{3})}, {\cal Z}; x_{1}, x_{2}),$$ 
$\alpha\in {\Bbb A}$, are nonzero, and
there exist 
$$F_{\alpha}(w'_{(a_{4})}, w_{(a_{1})},
w_{(a_{2})}, w_{(a_{3})}, {\cal Z}; x_{1}, x_{2})\in 
{\Bbb C}[x_{1}, x_{1}^{-1}, x_{2}, x_{2}^{-1}, (x_{1}-x_{2})^{-1}]$$
for $\alpha\in {\Bbb A}$,
such that  
\begin{eqnarray*}
\lefteqn{({\bf P}({\cal Z}))(w_{(a_{1})},
w_{(a_{2})}, w_{(a_{3})}; x_{1}, x_{2})}\nno\\
&&=\sum_{\alpha=1}^{m}g^{a_{1}, a_{2}, a_{3},
a_{4}}_{\alpha}(w_{(a_{1})},
w_{(a_{2})}, w_{(a_{3})}, {\cal Z}; x_{1}, x_{2})\iota_{12}\left(f^{a_{1}, 
a_{2}, a_{3}, a_{4}}_{\alpha}\right),
\end{eqnarray*}
and
\begin{eqnarray*}
\lefteqn{\langle w'_{(a_{4})}, g^{a_{1}, a_{2}, a_{3},
a_{4}}_{\alpha}(w_{(a_{1})},
w_{(a_{2})}, w_{(a_{3})}, {\cal Z}; x_{1}, x_{2})\rangle_{W_{a_{4}}}}\nno\\
&&=
\iota_{12}(F_{\alpha}(w'_{(a_{4})}, w_{(a_{1})},
w_{(a_{2})}, w_{(a_{3})}, {\cal Z}; x_{1}, x_{2}))
\end{eqnarray*} 
for $\alpha\in {\Bbb A}$.

{\bf (b) (commutativity in formal variables)}
We have
\begin{eqnarray*}
\lefteqn{\langle w'_{(a_{4})}, g^{a_{2}, a_{1}, a_{3},
a_{4}}_{\alpha}(w_{(a_{2})},
w_{(a_{1})}, w_{(a_{3})}, {\cal B}({\cal Z}); x_{2}, x_{1})
\rangle_{W_{a_{4}}}}\nno\\
&&=
\iota_{21}(F_{\alpha}(w'_{(a_{4})}, w_{(a_{1})},
w_{(a_{2})}, w_{(a_{3})}, {\cal Z}; x_{1}, x_{2}))
\end{eqnarray*}
for $\alpha\in {\Bbb A}$.
\end{theo}
\pf
By Lemma \ref{correl1}, we know that there exists $f\in {\Bbb G}^{a_{1},
a_{2}, a_{3}, a_{4}}$ such that  
$$\langle w'_{(a_{4})}, ({\bf P}({\cal Z}))(w_{(a_{1})},
w_{(a_{2})}, w_{(a_{3})}; x_{1}, x_{2})\rangle_{W_{a_{4}}}
=\iota_{12}f.$$
Since $\{f^{a_{1}, a_{2}, a_{3}, a_{4}}_{\alpha}\}_{\alpha\in
{\Bbb A}}$ is a basis of ${\Bbb G}^{a_{1},
a_{2}, a_{3}, a_{4}}$, there exist unique 
$$F_{\alpha}(w'_{(a_{4})}, w_{(a_{1})},
w_{(a_{2})}, w_{(a_{3})}, {\cal Z}; x_{1}, x_{2})\in {\Bbb
C}[x_{1}, x_{1}^{-1}, x_{2}, x_{2}^{-1}, (x_{1}-x_{2})^{-1}]$$
such that only finitely many of them are nonzero and 
$$f=\sum_{\alpha\in {\Bbb A}}F_{\alpha}(w'_{(a_{4})}, w_{(a_{1})},
w_{(a_{2})}, w_{(a_{3})}, {\cal Z}; x_{1}, x_{2})f^{a_{1}, 
a_{2}, a_{3}, a_{4}}_{\alpha}.$$
We define the maps $g^{a_{1},
a_{2}, a_{3}, a_{4}}_{\alpha}$ by
\begin{eqnarray*}
\lefteqn{\langle w'_{(a_{4})}, g^{a_{1},
a_{2}, a_{3}, a_{4}}_{\alpha}(w_{(a_{1})},
w_{(a_{2})}, w_{(a_{3})}, {\cal Z}; x_{1},
x_{2})\rangle_{W_{a_{4}}}}\nno\\ 
&&=F_{\alpha}(w'_{(a_{4})}, w_{(a_{1})},
w_{(a_{2})}, w_{(a_{3})}, {\cal Z}; x_{1}, x_{2})
\end{eqnarray*}
for $\alpha\in {\Bbb A}$. Then the generalized rationality for
products holds.

The conclusion {\bf (b)} follows immediately from the commutativity
for intertwining operators formulated in terms of complex variables.
\epfv

We also have the following generalized rationality of iterates 
and associativity
in formal variables:

\begin{theo}\label{associa}
There 
exist 
linear maps
\begin{eqnarray*}
\lefteqn{h^{a_{1}, a_{2}, a_{3},
a_{4}}_{\alpha}: W_{a_{1}}\otimes W_{a_{2}}\otimes W_{a_{3}}\otimes 
\coprod_{a_{5}\in {\cal A}}
{\cal V}_{a_{1}a_{2}}^{a_{5}}\otimes 
{\cal V}_{a_{5}a_{3}}^{a_{4}}}\nno\\
&&\hspace{8em} \to W_{a_{4}}[x_{0}, 
x_{0}^{-1}, x_{2}, x_{2}^{-1}][[x_{0}/x_{2}]]\nno\\
&&\hspace{3em} w_{(a_{1})}\otimes w_{(a_{2})}\otimes w_{(a_{3})}\otimes {\cal
Z}\nno\\
&&\hspace{8em}\mapsto h^{a_{1}, a_{2}, a_{3},
a_{4}}_{\alpha}(w_{(a_{1})},
w_{(a_{2})}, w_{(a_{3})}, {\cal Z}; x_{0}, x_{2}),
\end{eqnarray*}
$\alpha\in {\Bbb A}$,
satisfying the following property:  

{\bf (a) (generalized rationality of iterates)}
For any
$w_{(a_{1})}\in W_{a_{1}}$, $w_{(a_{2})}\in W_{a_{2}}$, $w_{(a_{3})}\in 
W_{a_{3}}$, $w'_{(a_{4})}\in W'_{a_{4}}$, and any
$${\cal Z}\in \coprod_{a_{5}\in {\cal A}}
{\cal V}_{a_{1}a_{2}}^{a_{5}}\otimes 
{\cal V}_{a_{5}a_{3}}^{a_{4}},$$
only finitely many of 
$$h^{a_{1}, a_{2}, a_{3},
a_{4}}_{\alpha}(w_{(a_{1})},
w_{(a_{2})}, w_{(a_{3})}, {\cal Z}; x_{0}, x_{2}),$$ 
$\alpha\in {\Bbb A}$, are nonzero, and
we have 
\begin{eqnarray*}
\lefteqn{
({\bf I}({\cal F}({\cal Z})))(w_{(a_{1})}, w_{(a_{2})}, 
w_{(a_{3})}; x_{0}, x_{2})}\nno\\
&&=\sum_{\alpha=1}^{m}h^{a_{1}, a_{2}, a_{3},
a_{4}}_{\alpha}(w_{(a_{1})},
w_{(a_{2})}, w_{(a_{3})}, {\cal F}({\cal Z}); x_{0}, x_{2})\iota_{20}
(f^{a_{1}, 
a_{2}, a_{3}, a_{4}}_{\alpha})
\end{eqnarray*}
and 
\begin{eqnarray*}
\lefteqn{\langle w'_{(a_{4})}, h^{a_{1}, a_{2}, a_{3},
a_{4}}_{\alpha}(w_{(a_{1})},
w_{(a_{2})}, w_{(a_{3})}, {\cal F}({\cal Z}); x_{0},
x_{2})\rangle_{W_{a_{4}}}}\nno\\
&&=\iota_{20}(F_{\alpha}(w'_{(a_{4})}, w_{(a_{1})},
w_{(a_{2})}, w_{(a_{3})}, {\cal Z}; x_{2}+x_{0}, x_{2}))
\end{eqnarray*}
for $\alpha\in {\Bbb A}$.

{\bf (b) (associativity in formal variables)} We have 
\begin{eqnarray*}
&\iota_{20}^{-1}\langle w'_{(a_{4})}, h^{a_{1}, a_{2}, a_{3},
a_{4}}_{\alpha}(w_{(a_{1})},
w_{(a_{2})}, w_{(a_{3})}, {\cal F}({\cal Z}); x_{0},
x_{2})\rangle_{W_{a_{4}}}\mbar_{x_{0}=x_{1}-x_{2}}&\nno\\
&=F_{\alpha}(w'_{(a_{4})}, w_{(a_{1})},
w_{(a_{2})}, w_{(a_{3})}, {\cal Z}; x_{1}, x_{2})&
\end{eqnarray*}
for $\alpha\in {\Bbb A}$.

\end{theo}
\pf
This result follows immediately from the associativity for
intertwining operators in terms of complex variables and Theorem 
\ref{rat}.
\epfv

\begin{rema}
{\rm These results are only about certain special ``four-point
functions'' (the four points being $\infty, z_{1}, z_{2}, 0$).  The
same method can be used to prove the generalized rationality property
of ``$n$-point functions'' and in particular the generalized
rationality property of general ``four-point functions.''  In the case
that the intertwining operator algebra is the one associated to an
affine Lie algebra and $w_{(a_{1})}, w_{(a_{2})}, w_{(a_{3})}$, and
$w'_{(a_{4})}$ are highest-weight vectors, it has been shown that 
the products and iterates of intertwining operators can be 
expressed using hypergeometric functions
(see, e.g., \cite{TK} and \cite{EFK}).}
\end{rema}

\begin{rema}\label{weaker-th}
{\rm Note that the proof of the results above, especially Lemma
\ref{correl1}, do not need the absolute convergence of products of
arbitrary numbers of intertwining operators in the corresponding
regions; the absolute convergence of
\begin{eqnarray}\label{conv-3sy}
\lefteqn{\langle w'_{(a_{4})}, {\cal Y}_{1}(w_{(a_{1})}, x_{1})
{\cal Y}_{2}(w_{(a_{2})}, x_{2})\cdot}\nno\\
&&\hspace{3em}\cdot e^{x_{3}L(-1)}w_{(a_{3})}\rangle_{M_{a_{4}}}
\lbar_{x_{1}^{n}= e^{n\log z_{1}}, \;x_{2}^{n}=e^{n\log z_{2}},
x_{3}=z_{3}}
\end{eqnarray}
in the region $|z_{1}|>|z_{2}|>|z_{3}|>0$ is enough. Thus for an
algebra
satisfying all the axioms for intertwining operator algebras except
that the absolute convergence of products of intertwining operators in
the corresponding regions is
weakened to the absolute  convergence of (\ref{conv-3sy}) in the
region $|z_{1}|>|z_{2}|>|z_{3}|>0$, the conclusions of 
Theorems \ref{rat} and \ref{associa} hold.} 
\end{rema}

One immediate consequence of Theorems \ref{rat} and Remark
\ref{weaker-th} 
is Proposition \ref{c&m->r}:

\vspace{1em}

\noindent {\it Proof of Proposition \ref{c&m->r}}: Note that an algebra 
satisfying the axioms for vertex operator algebras except for the Jacobi
identity, commutativity, associativity and the absolute convergence of
(\ref{conv-3y}) is an algebra discussed in Remark \ref{weaker-th}. 
Thus the conclusion of Theorems \ref{rat} and \ref{associa} hold. 
But  the vertex operators for such an algebra
 are  single-valued and so are their products 
and iterates. Thus the correlation functions must be rational.
\epfv

The Jacobi identity is also an easy consequence of the results we have
proved.

\vspace{1em}

\noindent {\it Proof of Theorem \ref{jacobi}}: The first two conclusions in 
Theorem \ref{jacobi} are parts of the conclusions of 
Theorems \ref{rat} and \ref{associa}. Using Theorems \ref{rat} and
\ref{associa} and Proposition 3.1.1 in \cite{FHL}, We obtain the
Jacobi identity. 
\epfv

In general, the multivalued functions corresponding to 
products and iterates of intertwining operators do not form a finitely
generated modules over ${\Bbb C}[x_{1}, x_{1}^{-1}, x_{2}, x_{2}^{-1},
(x_{1}-x_{2})^{-1}]$.  But when the intertwining operator algebra is
the intertwining operator algebra constructed from 
a vertex operator algebra 
generated by weight-one elements (i.e., a vertex operator
algebra associated to an
affine Lie algebra), these multivalued functions indeed form a finitely
generated modules. To prove this fact, we need:

\begin{lemma}
Assume that the intertwining operator algebra is
the intertwining operator algebra constructed from 
a vertex operator algebra $V$
generated by weight-one elements.
Let $\{{\cal Y}^{a_{3}}_{a_{1}a_{2};i}\}_{i=1}^{{\cal
N}_{a_{1}a_{2}}^{a_{3}}}$ be a basis of the space ${\cal
V}_{a_{1}a_{2}}^{a_{3}}$ for $a_{1}, a_{2}, a_{3}\in {\cal A}$,
$w_{(a), k}$, $k=1, \dots, \dim (W_{a})_{(n_{a})}$, and $w'_{(a), k}$,
$k=1, \dots, \dim (W_{a})_{(n_{a})}$ (note that 
$\dim (W'_{a})_{(n_{a})}=\dim (W_{a})_{(n_{a})}$), basis of
the lowest weight nonzero homogeneous subspace of $W_{a}$ and
$W'_{a}$, respectively, for $a\in {\cal A}$.  Then for any $a_{1},
a_{2}, a_{3}\in {\cal A}$, $w_{(a_{1})}\in W_{a_{1}}$, $w_{(a_{2})}\in
W_{a_{2}}$, $w_{(a_{3})}\in W_{a_{3}}$, $w'_{(a_{4})}\in W'_{a_{4}}$
and
$${\cal Z}\in \coprod_{a_{5}\in {\cal A}}
{\cal V}_{a_{1}a_{5}}^{a_{4}}\otimes 
{\cal V}_{a_{2}a_{3}}^{a_{5}},$$ 
$$\langle w'_{(a_{4})}, ({\bf P}({\cal Z}))(w_{(a_{1})},
w_{(a_{2})}, w_{(a_{3})}; x_{1}, x_{2})\rangle_{W_{a_{4}}}$$
is a linear combination of 
$$\langle w'_{(a_{4}), k_{4}}, {\cal
Y}^{a_{4}}_{a_{1}a_{5};i}(w_{(a_{1}), k_{1}}, x_{1}) 
{\cal
Y}^{a_{5}}_{a_{2}a_{3};i}(w_{(a_{2}), k_{2}}, x_{2})w_{(a_{3}),
k_{3}}\rangle_{W_{a_{4}}},$$ 
$a_{5}\in {\cal A}$, $i=1, \dots, {\cal
N}_{a_{1}a_{5}}^{a_{4}}$, $j=1, \dots, {\cal
N}_{a_{2}a_{3}}^{a_{5}}$, $k_{l}=1, \dots, \dim
(W_{a_{l}})_{(n_{a_{l}})}$ for $l=1, \dots, 4$, 
with coefficients in $\iota_{12}{\Bbb C}[x_{1}, x_{1}^{-1}, x_{2},
x_{2}^{-1}, (x_{1}-x_{2})^{-1}]$.
\end{lemma}
\pf
Since 
$$\langle w'_{(a_{4})}, ({\bf P}({\cal Z}))(w_{(a_{1})},
w_{(a_{2})}, w_{(a_{3})}; x_{1}, x_{2})\rangle_{W_{a_{4}}}$$
is a linear combination of 
$$\langle w'_{(a_{4})}, {\cal
Y}^{a_{4}}_{a_{1}a_{5};i}(w_{(a_{1})}, x_{1}) 
{\cal
Y}^{a_{5}}_{a_{2}a_{3};i}(w_{(a_{2})}, x_{2})w_{(a_{3})}
\rangle_{W_{a_{4}}}$$ 
with complex coefficients, we need only show that 
the lemma is true for such an expression.

Consider 
$$\langle w'_{(a_{4})}, {\cal
Y}^{a_{4}}_{a_{1}a_{5};i}(w_{(a_{1})}, x_{1}) 
{\cal
Y}^{a_{5}}_{a_{2}a_{3};i}(w_{(a_{2})}, x_{2})w_{(a_{3})}
\rangle_{W_{a_{4}}}.$$ 
Using the Jacobi identity defining intertwining operators, this
expression can be written as a linear combination of expressions of
the same form but with $w'_{(a_{4})}$, $w_{(a_{1})}$, $w_{(a_{2})}$
and $w_{(a_{3})}$ in the lowest weight spaces, and with elements of
$\iota_{12}{\Bbb C}[x_{1}, x_{1}^{-1}, x_{2}, x_{2}^{-1}, 
(x_{1}-x_{2})^{-1}]$
as coefficients. So the conclusion is true.
\epfv

The following result is an immediate consequence:

\begin{propo}
Assume that the intertwining operator algebra is
the intertwining operator algebra constructed from 
a vertex operator algebra $V$
generated by weight-one elements.
Then for $a_{1}, a_{2}, a_{3}, a_{4}\in {\cal A}$, 
there exist finitely many elements $f^{a_{1}, a_{2}, a_{3},
a_{4}}_{i}$, $i=1, \dots, m$, of ${\Bbb G}^{a_{1}, a_{2}, a_{3},
a_{4}}$, such that the space of multivalued functions obtained from 
products and iterates of intertwining operators is a free module over 
${\Bbb C}[x_{1}, x_{1}^{-1}, x_{2}, x_{2}^{-1}, 
(x_{1}-x_{2})^{-1}]$ with a basis $f^{a_{1}, a_{2}, a_{3},
a_{4}}_{i}$, $i=1, \dots, m$. \epf
\end{propo}

\renewcommand{\theequation}{\thesection.\arabic{equation}}
\renewcommand{\therema}{\thesection.\arabic{rema}}
\setcounter{equation}{0}
\setcounter{rema}{0}

\section{Intertwining operator algebras in terms of the Jacobi identity}

In this section, we give two equivalent 
definitions of intertwining operator algebra
using the Jacobi identity (\ref{1.0}) as the
main axiom. Intertwining operator algebras defined in this section are
natural nonabelian generalizations of abelian intertwining algebras 
introduced by Dong and Lepowsky (see \cite{DL1} and \cite{DL}).
The detailed axiomatic study of intertwining operator
algebras satisfying these
definitions, and in particular, the equivalence of these definitions
with each other and with
the one in
\cite{H4},  will be given in the forthcoming paper \cite{H6}.

\begin{defi}\label{4-1}
{\rm An {\it intertwining operator algebra of central charge 
$c\in {\Bbb C}$} is a vector space
$$W=\coprod_{a\in
{\cal A}}W^{a}$$
graded
by a finite set ${\cal A}$
containing a special element $e$
(graded  by {\it color}), 
equipped with a vertex operator algebra
structure of central charge $c$ 
on $W^{e}$, a $W^{e}$-module structure on $W^{a}$ for 
each $a\in {\cal A}$, a subspace ${\cal V}_{a_{1}a_{2}}^{a_{3}}$ of 
the space of all intertwining operators of type 
${W^{a_{3}}\choose W^{a_{1}}W^{a_{2}}}$ for  each triple
$a_{1}, a_{2}, a_{3}\in {\cal A}$, a {\it fusing isomorphism} 
$${\cal F}(a_{1}, a_{2}, a_{3}; a_{4}): 
\coprod_{a\in {\cal A}}{\cal V}_{a_{1}a}^{a_{4}}\otimes 
{\cal V}_{a_{2}a_{3}}^{a}
\to \coprod_{a\in {\cal A}}{\cal V}_{a_{1}a_{2}}^{a}\otimes {\cal
V}_{aa_{3}}^{a_{4}}$$ 
for each quadruple $a_{1}, a_{2}, a_{3}, a_{4}\in {\cal A}$,
 and a
{\it skew-symmetry isomorphism}  
$$\Omega(a_{1}, a_{2}; a_{3}): 
{\cal V}_{a_{1}a_{2}}^{a_{3}}\to {\cal
V}_{a_{2}a_{1}}^{a_{3}}$$
for each triple $a_{1}, a_{2}, a_{3}\in {\cal A}$,
satisfying the following axioms:

\begin{enumerate}

\item For any $a\in {\cal A}$, there exists $h_{a}\in {\Bbb R}$ 
such that the ${\Bbb C}$-graded module $W^{a}$
is $h_{a}+{\Bbb Z}$-graded.

\item The $W^{e}$-module structure on $W^{e}$ is the adjoint module
structure. For any $a\in {\cal A}$, the space ${\cal V}_{ea}^{a}$ is the
one-dimensional vector space spanned by the vertex operators for the
$W^{e}$-module $W^{a}$,  and for any ${\cal Y}\in {\cal V}_{ea}^{a}$
and any $w_{(e)}\in W^{e}$,
$w_{(a)}\in W^{a}$, 
$$((\Omega(e, a; a))({\cal Y}))(w_{(a)}, x)w_{(e)}=e^{xL(-1)}{\cal Y}(w_{(e)},
-x)w_{(a)}.$$ 
For any $a_{1}, a_{2}\in {\cal A}$ such that
$a_{1}\ne a_{2}$, ${\cal V}_{ea_{1}}^{a_{2}}=0$. 

\item  The genus-zero Moore-Seiberg equations:
Let 
$${\cal F}:
\coprod_{a_{1}, a_{2}, a_{3}, a_{4},
a_{5}\in {\cal A}}{\cal V}_{a_{1}a_{5}}^{a_{4}}\otimes 
{\cal V}_{a_{2}a_{3}}^{a_{5}}
\to \coprod_{a_{1}, a_{2}, a_{3}, a_{4},
a_{5}\in {\cal A}}{\cal V}_{a_{1}a_{2}}^{a_{5}}\otimes {\cal
V}_{a_{5}a_{3}}^{a_{4}}$$ 
and
$$\Omega: \coprod_{a_{1}, a_{2}, a_{3}\in {\cal A}}
{\cal V}_{a_{1}a_{2}}^{a_{3}}\to \coprod_{a_{1}, a_{2}, a_{3}\in {\cal A}}
{\cal V}_{a_{1}a_{2}}^{a_{3}}$$
be the isomorphisms obtained {}from 
${\cal F}(a_{1}, a_{2}, a_{3}; a_{4})$ and 
$\Omega(a_{1}, a_{2}; a_{3})$, $a_{1},$ $a_{2}, 
a_{3}, a_{4}\in {\cal A}$,
and let $F^{(1)}_{12}$, $F^{(2)}_{12}$, $F_{13}$, $F^{(1)}_{23}$,
$F^{(2)}_{23}$, $\Omega^{(p)}$, $(\Omega^{-1})^{(p)}$, $p=1, 2, 3, 4$, be
the isomorphisms obtained {}from ${\cal F}$ and $\Omega$ as in Section 3.
Then (\ref{pentagon}), (\ref{hexagon1}) and 
(\ref{hexagon2}) hold.

\item The Jacobi identity: Let 
${\cal B}(a_{1}, a_{2}; a_{3}, a_{4})$,
$a_{1}, a_{2}, a_{3}, a_{4}\in {\cal A}$,
be the braiding isomorphisms defined by (\ref{braiding}), 
$${\cal B}: \coprod_{a_{1}, a_{2}, a_{3}, a_{4}, a_{5}\in {\cal
A}}{\cal V}_{a_{1}a_{5}}^{a_{4}}\otimes {\cal V}_{a_{2}a_{3}}^{a_{5}}
\to \coprod_{a_{1}, a_{2}, a_{3}, a_{4}, a_{5}\in {\cal A}}{\cal
V}_{a_{1}a_{5}}^{a_{4}}\otimes {\cal V}_{a_{2}a_{3}}^{a_{5}}$$ the
braiding isomorphism obtained {}from ${\cal B}(a_{1}, a_{2}; a_{3},
a_{4})$, $a_{1}, a_{2}, a_{3}, a_{4}\in {\cal A}$, and
${\cal F}$ the isomorphism 
obtained above in the third axiom.
Then there exist linear maps $g^{a_{1}, a_{2}, a_{3}, a_{4}}_{\alpha}$
and $h^{a_{1}, a_{2}, a_{3}, a_{4}}_{\alpha}$ of the form 
(\ref{1.-4}) and (\ref{1.-3}) such that (\ref{1.-2}), (\ref{1.-1}) and
the the Jacobi identity
(\ref{1.0}) holds for any $a_{1}, a_{2}$, $a_{3}, a_{4}\in {\cal A}$, any
$w_{(a_{1})}\in W_{a_{1}}$, $w_{(a_{2})}\in W_{a_{2}}$, $w_{(a_{3})}\in 
W_{a_{3}}$, any
$${\cal Z}\in \coprod_{a_{4}, a_{5}\in {\cal A}}
{\cal V}_{a_{1}a_{5}}^{a_{4}}\otimes 
{\cal V}_{a_{2}a_{3}}^{a_{5}}\subset \coprod_{a_{1}, a_{2}, a_{3}, a_{4},
a_{5}\in {\cal A}}{\cal V}_{a_{1}a_{5}}^{a_{4}}\otimes 
{\cal V}_{a_{2}a_{3}}^{a_{5}},$$
and any $\alpha\in {\Bbb A}$.

\end{enumerate}}

\end{defi}

This definition is simple but it requires that 
the reader  be familiar with 
the notions of 
vertex operator algebra, module for a vertex operator algebra,
 and intertwining operator for a vertex operator algebra. 
The second definition is longer but much
more explicit, and it does not assume that the 
reader has any knowledge of the theory of vertex operator algebras.

Let $A$ be an $n$-dimensional commutative associative algebra over 
${\Bbb C}$. Then
for any basis ${\cal A}$ of $A$, there are structure constants ${\cal 
N}_{a_{1}a_{2}}^{a_{3}}\in {\Bbb C}$, 
$a_{1}, a_{2}, a_{3}\in {\cal A}$, such that 
$$a_{1}a_{2}=\sum_{a_{3}\in {\cal A}}{\cal 
N}_{a_{1}a_{2}}^{a_{3}}a_{3}$$ for any $a_{1}, a_{2}\in {\cal A}$.
Assume that $A$ has a basis ${\cal A}\subset A$ containing the identity
$e\in A$ such that all the structure constants ${\cal 
N}_{a_{1}a_{2}}^{a_{3}}$, $a_{1}, a_{2}, a_{3}\in {\cal A}$, are in
${\Bbb N}$.
Note that for any $a_{1}, a_{2}\in {\cal A}$, 
$${\cal N}_{ea_{1}}^{a_{2}}=\delta_{a_{1},
a_{2}}=\left\{\begin{array}{ll}
1&a_{1}=a_{2},\\0& a_{1}\ne a_{2}.\end{array}\right.$$
The commutativity and associativity of $A$ give the following
identities:
\begin{eqnarray*}
{\cal N}_{a_{1}a_{2}}^{a_{3}}&=&{\cal N}_{a_{2}a_{1}}^{a_{3}}\\
\sum_{a\in {\cal A}}{\cal N}_{a_{1}a_{2}}^{a}{\cal 
N}_{aa_{3}}^{a_{4}}&=&
\sum_{a\in {\cal A}}{\cal N}_{a_{1}a}^{a^{4}}{\cal 
N}_{a_{2}a_{3}}^{a},
\end{eqnarray*}
for $a_{1}, a_{2}, a_{3}, a_{4}\in {\cal A}$.

For a vector space
$W=\coprod_{a\in {\cal A}, n\in {\Bbb R}}W^{a}_{(n)}$ doubly graded
by ${\Bbb R}$ and  ${\cal A}$, let 
\begin{eqnarray*}
W_{(n)}&=&\coprod_{a\in {\cal A}}W^{a}_{(n)}\\
W^{a}&=&\coprod_{n\in {\Bbb R}}W^{a}_{(n)}.
\end{eqnarray*}
Then 
$$W=\coprod_{n\in {\Bbb R}}W_{(n)}=\coprod_{a\in {\cal A}}W^{a}.$$

\begin{defi}\label{4-2}
{\rm  An {\it intertwining operator algebra of central charge
$c\in {\Bbb C}$} consists of the following data:

\begin{enumerate}

\item A finite-dimensional commutative associative algebra $A$ and a
basis ${\cal A}$ of $A$ containing the identity $e\in A$ such that all
the structure constants ${\cal  N}_{a_{1}a_{2}}^{a_{3}}$, $a_{1},
a_{2}, a_{3}\in {\cal A}$, are in ${\Bbb N}$.

\item   A vector space
$$W=\coprod_{a\in
{\cal A}, n\in {\Bbb R}}W^{a}_{(n)},\; \mbox{\rm for}\; w\in W^{a}_{(n)},\;
n=\wt w,\; a=\clr w$$
doubly graded
by ${\Bbb R}$ and ${\cal A}$
 (graded by {\it weight} and by {\it color}, respectively).

\item For each triple $(a_{1}, a_{2}, a_{3})\in {\cal A}\times
{\cal A}\times {\cal A}$, 
 an ${\cal N}_{a_{1}a_{2}}^{a_{3}}$-dimensional vector space
${\cal V}_{a_{1}a_{2}}^{a_{3}}$ whose elements are linear maps
\begin{eqnarray*}
W^{a_{1}}&\to& \hom(W^{a_{2}}, W^{a_{3}})\{x\}\\
w_{(a_{1})}\in W^{a_{1}}&\mapsto &{\cal Y}(w_{(a_{1})}, x)\in 
\hom(W^{a_{2}}, W^{a_{3}})\{x\}.
\end{eqnarray*}

\item  For any $a_{1}, a_{2}, a_{3}, a_{4}\in {\cal A}$,
a {\it fusing isomorphism} 
$${\cal F}(a_{1}, a_{2}, a_{3}; a_{4}): 
\coprod_{a\in {\cal A}}{\cal V}_{a_{1}a}^{a_{4}}\otimes 
{\cal V}_{a_{2}a_{3}}^{a}
\to \coprod_{a\in {\cal A}}{\cal V}_{a_{1}a_{2}}^{a}\otimes {\cal
V}_{aa_{3}}^{a_{4}}$$ 
and for any $a_{1}, a_{2}, a_{3}\in {\cal A}$, a 
{\it skew-symmetry isomorphism}  
$$\Omega(a_{1}, a_{2}; a_{3}): 
{\cal V}_{a_{1}a_{2}}^{a_{3}}\to {\cal
V}_{a_{2}a_{1}}^{a_{3}}.$$

\item Two distinguished vectors ${\bf 1}\in W^{e}$ (the {\it  vacuum})
and $\omega\in W^{e}$ (the {\it Virasoro element}). 

\end{enumerate}

\noindent These data satisfy the
following axioms for $a_{1}, a_{2}, a_{3}\in {\cal A}$, 
$w_{(a_{1})}\in W^{a_{1}}$, and $w_{(a_{2})}\in W^{a_{2}}$:
\begin{enumerate}

\item The {\it weight-grading-restriction conditions}: 
For any $n\in {\Bbb Z}$ and $a\in {\cal A}$, 
$$\dim W^{a}_{(n)}< \infty$$
and for any $a\in {\cal A}$,
there exists $h_{a}\in {\Bbb R}$ such that 
$$W^{a}_{(n)}=0$$
for $n\not\in h_{a}+{\Bbb Z}$ or $n$
sufficiently small.

\item Axioms for intertwining operators:

\begin{enumerate}

\item The {\it single-valuedness condition}:  For any ${\cal Y}\in {\cal 
V}_{ea_{1}}^{a_{1}}$,
$${\cal Y}(w_{(a_{1})}, x)\in \mbox{Hom}(W^{a_{1}},
W^{a_{1}})[[x, x^{-1}]].$$

\item The {\it lower-truncation property for intertwining
operators}: For  any
${\cal Y}\in {\cal V}_{a_{1}a_{2}}^{a_{3}}$ and $w_{(a_{1})}\in W^{a_{1}}$, 
let ${\cal Y}_{n}(w_{(a_{1})})=\res_{x}x^{n}{\cal Y}(w_{(a_{1})}, x)$,
$n\in {\Bbb C}$, that is, ${\cal Y}(w_{(a_{1})}, x)=
\sum_{n\in {\Bbb R}}{\cal Y}_{n}(w_{(a_{1})})x^{-n-1}$. Then 
${\cal Y}_{n}(w_{(a_{1})})w_{(a_{2})}=0$ for $n$ sufficiently large.

\end{enumerate}

\item The genus-zero Moore-Seiberg equations (see the third axiom in
Definition \ref{4-1}).

\item Axioms for the vacuum:

\begin{enumerate}

\item The {\it identity property}: For any ${\cal Y}\in {\cal 
V}_{ea_{1}}^{a_{1}}$, there is $\lambda_{{\cal Y}}\in {\Bbb C}$ such that
${\cal Y}({\bf 1}, x)=\lambda_{{\cal Y}}I_{W^{a_{1}}}$, where
$I_{W^{a_{1}}}$ on the right is the identity
operator on $W^{a_{1}}$.

\item The {\it creation property}: For any ${\cal Y}\in {\cal
V}_{a_{1}e}^{a_{1}}$, there is $\mu_{{\cal Y}}\in {\Bbb C}$ such that ${\cal
Y}(w_{(a_{1})}, x){\bf 1}\in W[[x]]$ and $\lim_{x\to 0}{\cal Y}(w_{(a_{1})},
x){\bf 1}=\mu_{{\cal Y}}w_{(a_{1})}$ (that is, ${\cal Y}(w_{(a_{1})},
x){\bf 1}$ involves only nonnegative
integral powers of $x$ and the constant term is $\mu_{{\cal Y}}w_{(a_{1})}$).

\end{enumerate}

\item The Jacobi identity (see the fourth 
axiom in Definition \ref{4-1}).

\item Axioms for the Virasoro element:

\begin{enumerate}

\item The {\it Virasoro algebra relations}: Let $Y$ be the
element of ${\cal V}_{ea_{1}}^{a_{1}}$ such that $Y({\bf 1},
x)=I_{W^{a_{1}}}$ and let $Y(\omega, x)=\sum_{n\in
{\Bbb Z}}L(n)x^{-n-2}$. Then
$$[L(m), L(n)]=(m-n)L(m+n)+\frac{c}{12}(m^{3}-m)\delta_{m+n, 0}$$
for $m, n\in {\Bbb Z}$.

\item The {\it $L(0)$-grading property}: 
$L(0)w_{(a_{1})}=nw_{(a_{1})}=(\wt w_{(a_{1})})w_{(a_{1})}$ 
for $n\in {\Bbb R}$
and $w_{(a_{1})}\in W^{a_{1}}_{(n)}$.

\item The {\it $L(-1)$-derivative property}: For any ${\cal Y}\in {\cal 
V}_{a_{1}a_{2}}^{a_{3}}$,
$$\frac{d}{dx}{\cal Y}(w_{(a_{1})}, x)={\cal Y}(L(-1)w_{(a_{1})}, x).$$

\end{enumerate}

\end{enumerate}}
\end{defi}

{\small \sc Department of Mathematics, Rutgers University,
110 Frelinghuysen Rd., Piscataway, NJ 08854-8019}

{\em E-mail address}: yzhuang@math.rutgers.edu

\end{document}